\documentclass[a4paper, english, 11pt, onecolumn, oneside,final]{article}	
\usepackage[verbose, a4paper, hcentering]{geometry}	

\usepackage[utf8]{inputenc} 
\usepackage[T1]{fontenc}    
\usepackage[pdftitle={Sinkhorn}, pdfauthor={Alois Pichler}, citecolor=gray, urlcolor= gray, linkcolor= gray, unicode=true,  bookmarks=true, bookmarksnumbered=false, bookmarksopen=false, breaklinks=false, pdfborder={0 0 0}, pdfborderstyle={}, backref=page, colorlinks=true]{hyperref}
\usepackage[numbers]{natbib}	
\usepackage{booktabs}       
\usepackage{amsthm, mathtools, amssymb, nicefrac}
\mathtoolsset{showonlyrefs}
\DeclareMathOperator{\diag}{diag}

\usepackage{microtype}      
\usepackage{newtxtext, newtxmath}
\usepackage{dsfont}   
\usepackage[notcite,notref, draft]{showkeys}	

\usepackage[inline, shortlabels]{enumitem}           
\setlist[enumerate,1]{label=(\roman*)}  
\newcommand{\one}{{\mathbf 1}} 		
\newcommand{\zero}{{\mathbf 0}} 		
\usepackage{xcolor}         
\usepackage{siunitx}
\usepackage{todonotes}
\usepackage[export]{adjustbox}

\usepackage{subfig}
\usepackage{graphicx}

\usepackage{pgfplots}
\pgfplotsset{compat=newest}
\usepackage{siunitx}
\usepackage{tikz} 
\usetikzlibrary{trees,matrix,arrows.meta}
\usepackage[ruled, vlined]{algorithm2e}			
\usepackage{orcidlink}
\usepackage{diagbox}
\usepackage{doi}
\usepackage{intcalc}
\usepackage{rerunfilecheck}
\usepackage{footnote}
\makesavenoteenv{tabular}
\makesavenoteenv{table}

\newcounter{relctr} 
\everydisplay\expandafter{\the\everydisplay\setcounter{relctr}{0}} 

\AtBeginDocument{}

\theoremstyle{plain}
	\newtheorem{theorem}{Theorem}[section]		
		
	\newtheorem{lemma}[theorem]{Lemma}
	\newtheorem{proposition}[theorem]{Proposition}
\theoremstyle{definition}
	\newtheorem{definition}[theorem]{Definition}
\theoremstyle{remark}
	\newtheorem{remark}[theorem]{Remark}

\title{Nonequispaced Fast Fourier Transform Boost for the Sinkhorn~Algorithm}

\author{%
	Rajmadan Lakshmanan\thanks{Faculty of Mathematics, University of Technology, Chemnitz, Germany}\,
	\footnote{\href{rajmadan.lakshmanan@math.tu-chemnitz.de}{rajmadan.lakshmanan@math.tu-chemnitz.de}}
	\and
	Alois~Pichler\footnotemark[1]\, 
	\thanks{DFG, German Research Foundation – Project-ID 416228727 – SFB 1410.}\,
	\footnote{
		\href{https://orcid.org/0000-0001-8876-2429}{\orcidlink{0000-0001-8876-2429} https://orcid.org/0000-0001-8876-2429}}
	\and
	Daniel~Potts\footnotemark[1]\, \footnotemark[3]\,
	\thanks{
	\href{https://orcid.org/0000-0003-3651-4364} {\orcidlink{0000-0003-3651-4364} https://orcid.org/0000-0003-3651-4364}
}}

\begin{document}

\maketitle
\begin{abstract}
	This contribution features an accelerated computation of the Sinkhorn's algorithm, which approximates the Wasserstein transportation distance, by employing nonequispaced fast Fourier transforms (NFFT). 
	The algorithm proposed allows approximations of the Wasserstein distance by involving not more than $\mathcal O(n\log n)$ operations for probability measures supported by~$n$ points.
	Furthermore, the proposed method avoids expensive allocations of the characterizing matrices.
	With this numerical acceleration, the transportation distance is accessible to probability measures out of reach so far.
	Numerical experiments using synthetic and real data affirm the computational advantage and superiority.
	\medskip


\end{abstract}

\section{Introduction}\label{sec:intro}


In optimal transport theory, the Wasserstein distance~-- often referred to as the Monge-Kantorovich distance~-- is used to define and quantify optimal transitions between probability measures.
The lowest (or cheapest) average costs to fully transfer one probability measure into another characterizes the distance.
In most applications, costs correspond to the distance between locations.
For a comprehensive discussion of the Wasserstein distance from mathematical perspective we may refer to 
\citet{Villani2009}.

The concept of entropy regularization of the Wasserstein distance, proposed by \citet{cuturi2013sinkhorn}, is an important touchstone, which improves the computational process of traditional methods. 
This entropy regularized Wasserstein problem is  efficiently solved using the Sinkhorn's algorithm (cf.\ \citet{Sinkhorn1967a}).
In today's data-driven world, the powerful and growing relationship between optimization and data science utilizes the Wasserstein distance, e.g., for text classification (cf.\ \citet{kusner2015word}), clustering (cf.\ \citet{chakraborty2020hierarchical}), image classification (cf.\ \citet{pmlr-v139-tai21a}) or domain adaptation (cf.\ \citet{courty2014domain}).
Notably, most of the applications rely on \emph{discrete} measures. 
However, some significant contributions are also presented in literature to support the arguments of semi-discrete and/\,or continuous measures (cf.\ \citet{mensch2020online}). 
The constructive line of research on the entropy regularization method to approximate the Wasserstein distance proposes many significant algorithms to increase the computational efficiency, as well as to stabilize the approximation accuracy (cf.\ \citet{dvurechensky2018computational}, \citet{lin2019efficiency} or \citet{schmitzer2019stabilized}). 
However, this article addresses the efficient computation of standard Sinkhorn's algorithm in terms of time and memory allocation to approximate the  Wasserstein distance in a simple personal computer, especially in case of large data volume.
\paragraph{Related works.}
In the past decade, based on the well-known standard (equispaced) fast Fourier transform~(FFT) method, many approaches have been proposed to find efficient data representation for various problems. 
The family of standard FFT algorithms has been applied in many areas, such as face recognition (cf.\ \citet{hao2013facial}), autonomous vehicles (cf.\ \citet{bilik2019rise}), voice assistants (cf.\ \citet{revay2019multiclass}), etc., and have achieved notable performance. 
The standard FFT algorithm improves the computational operations from $\mathcal O(n^2)$  to $\mathcal O(n \log n)$, where $n$ denotes the number of data points, this process involves equispaced sampling. 
However, in some cases, the equispaced sampling  is one of the root causes of failure to meet accuracy (cf.\ \citet{platte2011impossibility}, \citet[Chapter~7]{PlPoStTa18}).   
We recognize that the optimal transport (OT) communities use the idea of standard FFT to speed up the Sinkhorn's iterations in some places (cf.\ \citet{papadakis2014optimal}). 
The standard FFT methods utilize equispaced convolution, which is a setback, when we consider the stability of the computation and approximation accuracy (cf.\ \citet[Remark~4.18]{peyre2019computational}).
To overcome this challenge we present a \emph{non}-equispaced convolution below, and it is also achievable in $\mathcal O(n \log n)$ arithmetic operations.
Furthermore, for faster computation, low-rank factorization techniques are considered as a popular argument among OT communities (cf.\ \citet{NEURIPS2019_f55cadb9,NEURIPS2020_9bde76f2,altschuler2020polynomial}). 
As a consequence of the line of research on Low-Rank Factorization, \citet{pmlr-v139-scetbon21a} have developed an algorithm to efficiently solve the regularized OT problem, which depends on \emph{low-rank couplings}.
As well, the method can be employed to accelerate problems involving multi marginals, cf.\ \citet{QuellmalzFFT}.

\paragraph{Contribution.}
We improve the computational time and memory allocation of the standard entropy regularization approach to approximate Wasserstein distance with negligible or no compromise of accuracy.
The technique we present here is a fast summation method, and it is based on the nonequispaced fast Fourier transform (NFFT), see \citet[Chapter~7]{PlPoStTa18}.
Using NFFT, we boost the performance of standard entropy regularization of Wasserstein distance with stable computation and high (machine) accuracy.
Additionally, we explicitly provide the bounds for the approximation of the Wasserstein distance.
We experimentally substantiate the computational efficiency of our proposed algorithm,  and we validate the accuracy via numerical results.

\paragraph{Outline of the paper.} This paper is organized as follows. 
Initially, in Section~\ref{sec:Preliminaries}, we discuss the necessary notations and definitions of Wasserstein distance. 
Section~\ref{sec:Sinkhorn divergence}  introduces the entropy regularization approach to approximate the Wasserstein distance~(Primal problem) and its dual formulation. 
Additionally, we show the convergence properties of Sinkhorn's iteration and recall the Sinkhorn divergence. A
fast summation technique based on NFFT, which is utilized
in this paper, is introduced in Section~\ref{sec:Fast_summation}. 
In Section~\ref{Sec:FFT_boost}, we propose the NFFT-accelerated Sinkhorn's algorithm and schematically explain the operations. 
Section~\ref{Sec:Numerical_exposition} contains the demonstration of performance of our proposed algorithm on synthetic as well as real data sets. 
Finally, Section~\ref{Sec:summary} summarizes and concludes the paper.   

\section{Preliminaries} \label{sec:Preliminaries}
In this section, we provide a short review of the  Monge–Kantorovich or the Wasserstein distance. 

On a space of probability measures, Wasserstein distances offer a natural metric. Intuitively, the Wasserstein distance measures the minimum, average amount of transporting cost required to transform one distribution into another. 

\begin{definition}[Wasserstein distance]\label{def:Wasserstein}
	Let $(\mathcal X,d)$ be a Polish space and $P$ and $\tilde P\in\mathcal P(\mathcal X) $ be two probability measures on the Borel sets of~$\mathcal X$. 
	The Wasserstein distance of order $r\ge1$ of the probability measures~$P$ and~$\tilde P $ for a given cost or distance function  $d \colon\mathcal X\times \mathcal X\to\mathbb R$  is
	\begin{align}
		W_r(P,\tilde P)&\coloneqq w_r(P,\tilde P)^{\nicefrac 1r},
	\shortintertext{where}
		w_r(P,\tilde P)&\coloneq \label{eq:W0}
		\inf_{\pi \in \Pi(P, \tilde P)} \iint_{\mathcal X\times{\mathcal X}} d(x,\,\tilde x)^r\, \pi(\mathrm d x,\mathrm d\tilde x).
	\end{align}
	Here, $\Pi(P,\tilde{P})\subset \mathcal P(\mathcal X^2)$ is the set of bivariate probability measures on $\mathcal X\times {\mathcal X}$ with marginals~$P$ and~$\tilde P$, respectively;
	that is, $\pi (A\times \mathcal X)=P(A)$ and $\pi (\mathcal X\times B)= \tilde P(B)$ for all measureables sets $A$ and $B\subset\mathcal X$.
\end{definition}
Wasserstein distances metrize the weak* topology on measures with finite $r$th moment. In the discrete setting considered below and all regular situations, the infimum in~\eqref{eq:W0} is attained (cf.\ \citet[Theorem~7.12]{Villani2003}).
\paragraph{Discrete framework.}

Concrete implementations of the Wasserstein problem rely on discrete measures of the form%
	\footnote{\(\delta_x(A)\coloneqq \one_A(x)=
	\begin{cases}
		1 & \text{if }x\in A,\\
		0 & \text{else}
	\end{cases}\) is the Dirac measure located at $x\in \mathcal X$.}
\begin{equation}\label{eq:Discrete}
	P(\cdot)= \sum_{i=1}^n p_i\, \delta_{x_i}(\cdot)
\end{equation}
with $p_i\ge 0$ and $\sum_{i=1}^n p_i=1$.
These measures are dense in $\mathcal P(\mathcal X)$ with respect to the weak* topology, see \citet{Bolley2008}.

For two discrete probability measures
\begin{equation}\label{eq:Prob}	P=\sum_{i=1}^np_i \delta_{x_i}\ \text{  and }\
	\tilde P=\sum_{j=1}^{\tilde n}\tilde p_j\,\delta_{\tilde x_j},
\end{equation}
the bivariate measure $\pi= \sum_{i=1}^n\sum_{j=1}^{\tilde n} \delta_{(x_i,\tilde x_j)}\in \mathcal P(\mathcal X^2)$ solves the Wasserstein problem~\eqref{eq:W0}, provided that the matrix $\pi=(\pi_{ij}) \in \mathbb R^{n\times\tilde n}$ is the solution of the optimization problem
	\begin{subequations}
		\begin{align}\label{eq:Da}
			w_r(P,\tilde P)=
			&\min
			\sum_{i=1}^n\sum_{j=1}^{\tilde n} \pi_{ij}\, d_{ij}^r,
		\intertext{where}
			&\sum_{j=1}^{\tilde n} \pi_{ij}= p_i  \text{ for }i= 1,\dots,n, \label{eq:Db}\\
			&\sum_{i=1}^n\pi_{ij}=\tilde p_j \text{ for }j= 1,\dots,\tilde n\ \text{ and} \label{eq:Dc}\\
			&\pi_{ij}\ge0 \text{ for all }i=1,\dots, n\text{ and } j=1,\dots,\tilde n	\label{eq:Dd}
		\end{align}
	\end{subequations}
	and $d_{ij}= d(x_i,x_j)$ is the distance matrix.
	The problem~\eqref{eq:Da}--\eqref{eq:Dd} is a linear optimization problem, occasionally referred to as \emph{Kantorovich problem}.
	In what follows, the optimal matrix is denoted~$\pi^w$.

\paragraph{Complexity.} For $n\approx \tilde n$, the linear optimization problem~\eqref{eq:Da}--\eqref{eq:Dd} can be solved by straightforward computation involving $\mathcal O(n^3)$ multiplications.

In the following Section~\ref{sec:Sinkhorn divergence} we recall the popular approach based on entropy regularization to reduce the computational burden of the optimization problem~\eqref{eq:Da}--\eqref{eq:Dd}.
This approach is efficiently tackled by an iteration process, which is popularly known as \emph{Sinkhorn's algorithm}. This algorithm is also known as matrix scaling type algorithm (cf.\ \citet{RoteZachariasen}). 
\paragraph{Notations.} Throughout this article, $\|\cdot\|$ stands for Euclidean norm or $2$-norm and $\|\cdot\|_1$ stands for $1$-norm. The vector of all ones and zeros denote as  $\one_n \coloneqq (1,\dots,1)^{\top} \in\mathbb R^n$ and $\zero_{n} \coloneqq (0,\dots,0)^{\top} \in\mathbb R^{{n}}$. For any probability vectors~$a$ and~$b$, the Kullback--Leibler divergence is
\[ D_\textit{KL}(a\mid b)\coloneqq \sum^n_{i=1} a_i \log\Big(\frac{a_i}{b_i}\Big).  \]  
\section{Entropy regularization and Sinkhorn divergences} \label{sec:Sinkhorn divergence}

This section considers the entropy-regularization of the Wasserstein problem, and characterizes its duality.
Furthermore, we recall Sinkhorn’s algorithm which permits a considerably faster implementation.

\subsection{Entropy regularization of the Wasserstein problem}
The Entropy Regularized Wasserstein~(ERW) problem involves the entropic regularization term
\begin{equation}\label{eq:negentropy}
	H(\pi)\coloneqq -\sum_{i,j}\pi_{ij}\log \pi_{ij}
\end{equation} to the linear optimization problem~\eqref{eq:Da}--\eqref{eq:Dd}.

\begin{definition}[ERW distance]\label{DefSinkhornDistance}
	The ERW distance with regularization parameter $\lambda>0$ is the minimal value of the optimization problem 
	\begin{subequations}
		\begin{align}\label{eq:Sinkhorn}
			s_{r;\lambda}(P,\tilde P)\coloneqq
			\min
			&\sum_{i=1}^n\sum_{j=1}^{\tilde n} \pi_{ij}\, d_{ij}^r- \frac1\lambda H(\pi),
		\shortintertext{where}
			& \sum_{j=1}^{\tilde n}\pi_{ij}= p_i  \text{ for }i= 1,\dots,n, \label{eq:SD1}\\
			& \sum_{i=1}^n\pi_{ij}=\tilde p_j \text{ for }j= 1,\dots,\tilde n\ \text{ and} \label{eq:SD2}\\
			& \pi_{ij}>0 \text{ for all }i= 1,\dots, n\text{ and } j= 1,\dots,\tilde n.	\label{eq:Sinkhornb}
		\end{align}
	\end{subequations}
	The matrix minimizing~\eqref{eq:Sinkhorn} subject to the constraintes \eqref{eq:SD1}--\eqref{eq:Sinkhornb} is denoted $\pi^{s} \in \mathbb{R}^{n\times\tilde n}$.
	Further, we set
	\begin{equation}\label{eq:15}
		\tilde s_{r;\lambda}\coloneqq  s_{r;\lambda}+\frac1\lambda H(\pi^s)= \sum_{i=1}^n\sum_{j=1}^{\tilde n} \pi_{ij}^s\,d_{ij}^r.
	\end{equation}
\end{definition}
The non-negativity constraint~\eqref{eq:Dd} is notably not active in the constraints~\eqref{eq:SD1}--\eqref{eq:Sinkhornb}, as the function $\varphi(x)\coloneqq x\log x$ is strictly convex in~$[0,1]$ with $\varphi^\prime(0)=-\infty$, and the optimal solution consequently satisfies $\pi_{ij}> 0$.
The regularizing term $\frac1\lambda H(\cdot)$ is strictly convex, so that the solution of the problem~\eqref{eq:Sinkhorn}--\eqref{eq:Sinkhornb} exists and is unique.
\begin{remark}[Regularizing term]
	To surpass the difficulty of numerical computation of the linear optimization problem~\eqref{eq:Da}--\eqref{eq:Dd}, the entropy regularization approach was originally proposed in \citet{cuturi2013sinkhorn}.
We also refer to \citet{gasnikov2016efficient}, which comprises the argument of efficient numerical methods for entropy linear programming problems.
\end{remark}

\paragraph{Choice of the regularization parameter $\lambda$.} 
In general, the selection of the regularization parameter~$\lambda$ plays a crucial role to obtain a good approximation of the Wasserstein distance.
From~\eqref{eq:Sinkhorn} and the arguments below, we infer that if $\lambda\to \infty$, we obtain the standard Wasserstein distance in the limit. 
We refer to \citet{neumayer2021optimal}, who study the regularization parameter.
The constructive line of research by \citet{feydy2020geometric} affirms that when the regularization parameter~$\lambda$ is not sufficiently large, the transportation plan and the regularized Wasserstein distance may be inconsiderable.
However, from the literature, we infer that the choice $\lambda\geq 20$ is a good bargain between accuracy and computational speed (cf.\ \citet[Figure~1.2, Figures~3.1--3.3]{genevay2019entropy}, \citet[Figure~2]{pmlr-v139-scetbon21a}, \citet[Section~7]{neumayer2021optimal}).\footnote[2]{ Cf.\ \href{https://marcocuturi.net/SI.html}{https://marcocuturi.net/SI.html}} 
To acquire a better approximation accuracy, we can increase the regularization parameter~$\lambda> 20$ with the price of relatively more arithmetic operations.
\medskip

In some applications it is crucial to estimate the Wasserstein distance with given accuracy.
This can be accomplished by choosing the regularization parameter~$\lambda$ large enough.
The following Lemma~\ref{lem:Sinkhorn Approximation} gives a precise instruction how to choose~$\lambda$ to obtain a prescribed accuracy. 
\begin{lemma}[Quality of the Sinkhorn Approximation]\label{lem:Sinkhorn Approximation}
	For $\varepsilon>0$ it holds that
	\begin{equation}
		s_{r;\lambda}(P,\tilde P)\le w_r(P,\tilde P)^r \le \tilde{s}_{r;\lambda}(P,\tilde P)\le s_{r;\lambda}(P,\tilde P)+ \varepsilon,	\label{Eq:SinkhornInequality}
	\end{equation}
	provided that
	\[	\lambda\ge \frac{H(P)+H(\tilde P)}{\varepsilon}.\]
	Here, $H(P)=-\sum_{i=1}^np_i\log p_i$ ($H(\tilde P)=-\sum_{j=1}^{\tilde n}\tilde p_j\log\tilde p_j$, resp.)\ is the entropy of the measure $P$ ($\tilde P$, resp.).
	Further, the entropies are bounded by $H(P)+H(\tilde P)\le \log n+\log\tilde n$.
\end{lemma}
\begin{proof}
	The first inequality in~\eqref{Eq:SinkhornInequality} follows by substituting the matrix $\pi^w$ in~\eqref{eq:Sinkhorn}, as $H(\pi^w)>0$.
	Further, with the matrix $\pi^s$, it holds that $w_{r;\lambda}(P,\tilde P)\le \tilde s_{r;\lambda}(P,\tilde P)$, which is the second inequality using~\eqref{eq:15}.

	Now let~$\pi$ be any matrix with marginals~$p$ (cf.~\eqref{eq:Db}) and~$\tilde p$ (cf.~\eqref{eq:Dc}).
	It follows with the log sum inequality (or Gibbs' inequality)
	\[	\sum_{i,j}\pi_{ij}\log\frac{\pi_{ij}}{p_i\,\tilde p_j}\ge 0\]
	that
	\begin{align}
		\sum_{i,j}\pi_{ij}\log\pi_{ij}&\ge\sum_{i,j}\pi_{ij}\log p_i +\sum_{i,j}\pi_{ij}\log\tilde p_j
		= \sum_{i}p_i\log p_i +\sum_j\tilde p_j\log\tilde p_j \\
		&= \sum_{i,j}p_i \tilde p_j\log p_i +\sum_{i,j}p_i\tilde p_j\log\tilde p_j = \sum_{i,j}p_i \tilde p_j \log(p_i\tilde p_j)\\
		&=\sum_i p_i \log p_i+\sum_j\tilde p_j\log \tilde p_j,
	\end{align}
	that is, $H(\pi)\le H(P)+H(\tilde P)$ and thus
	$\frac1\lambda H(\pi)\le \frac{H(P)+H(\tilde P)}\lambda\le \varepsilon$ for the parameter $\lambda$ large enough as in the assumption.

	The remaining inequality follows from $\tilde s_{r;\lambda}(P,\tilde P)= s_{r;\lambda}(P,\tilde P)+\frac1\lambda H(\pi^s)\le s_{r;\lambda}(P,\tilde P)+\varepsilon$.

	The inequality $H(P)\le\log n$ follows by applying Gibb’s inequality to the measures with weights~$p$ ($\tilde p$, resp.)\ and the constant weights $\big(\nicefrac1n,\dots,\nicefrac1n\big)$ ($\big(\nicefrac1{\tilde n},\dots,\nicefrac1{\tilde n}\big)$, resp.).
\end{proof}
\begin{remark}
	Note, that the choice $\lambda \ge \frac{\log n+\log\tilde n}\varepsilon$ is independent of the probability measure, but only depends on their granularity or dimension~$n$ ($\tilde n$, resp.).
	We may also refer to \citet[Proposition~1]{NEURIPS2018_3fc2c60b} and references therein for further, related inequalities for continuous measures.
\end{remark}

\subsubsection{Dual representation of entropy-regularized Wasserstein distance}
We restate the optimization problem of ERW in the following dual formulation.
\begin{proposition}[{cf.\ \citet[Remark~4.28]{peyre2019computational}}]
	For $\lambda>0$, the ERW~\eqref{eq:Sinkhorn}, \eqref{eq:SD1}--\eqref{eq:SD2} admits the following dual representation
	~\label{Prop:CovergenceofSinkhornandduality}
	\begin{equation}\label{eq:14}
		d(\alpha,\tilde\alpha)\coloneqq
		\max_{\alpha\in\mathbb R^n,\ \tilde\alpha\in\mathbb R^{\tilde n}}
	\frac1\lambda+ \frac1\lambda\sum_{i=1}^n p_i \log\alpha_i +\frac1\lambda\sum_{j=1}^{\tilde n}\tilde p_j\log\tilde\alpha_j-\frac1\lambda\sum_{i=1}^n\sum_{j=1}^{\tilde n} \alpha_i\, \mathrm{e}^{-\lambda\, d_{ij}^r}\, \tilde\alpha_j,
	\end{equation}
	which is strictly concave dual function.
\end{proposition}
\begin{proof}
	The Lagrangian of the ERW problem~\eqref{eq:Sinkhorn} with dual parameters $\beta$ (for the constraint~\eqref{eq:Db}) and $\gamma$ (for~\eqref{eq:Dc}) is
	\begin{align}\label{eq:12}
		\MoveEqLeft[6] L(\pi;\beta,\gamma)= \sum_{i=1}^n\sum_{j=1}^{\tilde n} d_{ij}^r\,\pi_{ij}+ \frac1\lambda \sum_{i=1}^n\sum_{j=1}^{\tilde n} \pi_{ij}\log \pi_{ij}
		+ \sum_{i=1}^n \beta_i\Big(p_i-\sum_{j=1}^{\tilde n} \pi_{ij}\Big)
		+ \sum_{j=1}^{\tilde n} \gamma_j\Big(\tilde p_j-\sum_{i=1}^n \pi_{ij}\Big).
	\end{align}
	The optimal measure~$\pi^*$ satisfying the first order constraint
	\begin{equation}
		0= \frac{\partial L}{\partial \pi_{ij}}= d_{ij}^r+ \frac1\lambda(\log \pi_{ij}+1)-\beta_i-\gamma_j 
	\end{equation}
	is
	\begin{equation}\label{eq:13}
		\pi_{ij}^*= \exp\big(-\lambda(d_{ij}^r- \beta_i- \gamma_j)- 1\big).
	\end{equation}
	The measure~$\pi^*$ minimizes the Lagrangian $L$ for $\beta$ and $\gamma$ fixed, and reveals the dual function
	\begin{align}
		d(\beta,\gamma)&=\inf_\pi L(\pi;\beta,\gamma)= L(\pi^*;\beta,\gamma)\\
		&= \sum_{i=1}^n\sum_{j=1}^{\tilde n} d_{ij}^r\, \pi_{ij}^*+ \frac1\lambda\sum_{i=1}^n\sum_{j=1}^{\tilde n} \pi_{ij}^*\big(-\lambda(d_{ij}^r-\beta_i-\gamma_j)-1\big)& \\
		&\qquad+ \sum_{i=1}^n\beta_i\Big(p_i-\sum_{j=1}^{\tilde n} \pi_{ij}^*\Big)
		+ \sum_{j=1}^{\tilde n}\gamma_j\Big(\tilde p_j-\sum_{i=1}^n\pi_{ij}^*\Big)\\
		&= -\frac1\lambda\sum_{i=1}^n\sum_{j=1}^{\tilde n} \pi_{ij}^*+ \sum_{i=1}^n\beta_i\,p_i	+ \sum_{j=1}^{\tilde n}\gamma_j\,\tilde p_j\\
		&= \sum_{i=1}^np_i\,\beta_i	+ \sum_{j=1}^{\tilde n}\tilde p_j\,\gamma_j -\frac1\lambda\sum_{i=1}^n\sum_{j=1}^{\tilde n} \mathrm{e}^{-\lambda(d_{ij}^r-\beta_i-\gamma_j)-1}\label{eq:11}
	\end{align}
	explicitly. Now substitute $\alpha_i= \mathrm{e}^{\lambda\,\beta_i-\nicefrac12}$ and $\tilde\alpha_j= \mathrm{e}^{\lambda\,\gamma_j-\nicefrac12}$,
	then the dual function is 
	\begin{equation}\label{eq:27}
		d(\alpha,\tilde\alpha)= \frac1\lambda\sum_{i=1}^n p_i\, \Big(\frac12+\log\alpha_i\Big) +\frac1\lambda\sum_{j=1}^{\tilde n}\tilde p_j\,\Big(\frac12+\log\tilde\alpha_j\Big)-\frac1\lambda\sum_{i=1}^n\sum_{j=1}^{\tilde n} \alpha_i\, \mathrm{e}^{-\lambda\, d_{ij}^r}\, \tilde\alpha_j.
	\end{equation}
	The assertion of the proposition thus follows, as $\sum_{i=1}^n p_i=1$ and  $\sum_{j=1}^{\tilde n} \tilde p_j=1$ and as the duality gap vanishes for the strictly convex objective function~\eqref{eq:Sinkhorn}.
\end{proof}
\paragraph{Scaling variables and kernel matrix.}
Indeed, first of all, we notice that optimal measure~\eqref{eq:13} can be obtained in terms of the scaling variables $\alpha_i$ and $\tilde\alpha_j$ by 
\begin{equation}\label{eq:pi*}
	\pi^*= \diag(\alpha)\, \mathrm{e}^{-\lambda\, d^r} \, \diag(\tilde\alpha).
\end{equation}
The aforementioned dual problem~\eqref{eq:14} can be solved by a matrix scaling algorithm, which is popularly known as Sinkhorn's algorithm. 
Further, the derivate of~\eqref{eq:14} with respect to $\alpha_i$ ($\tilde\alpha_i$, resp.)\ gives the first order conditions, which is the basis for Sinkhorn's iteration and it is expressed as 
\begin{equation}\label{eq:S1}
	\alpha_i\coloneqq \frac{p_i}{\sum_{j=1}^{\tilde n}  \mathrm{e}^{-\lambda\,d_{ij}^r}\,\tilde\alpha_j} \text{ and } \tilde\alpha_j\coloneqq \frac{p_j}{\sum_{i=1}^{n}\,\mathrm{e}^{-\lambda\,d_{ij}^r}\,\alpha_i}.
\end{equation}
The main computational bottleneck of Sinkhorn's iterations is the matrix-vector multiplication in~\eqref{eq:S1}, which requires $\cal O$$(n\cdot\tilde n)$ arithmetic operations.
In our study, we relax the computational burden by taking advantage of special structure of the matrix \[k_{ij}\coloneqq \mathrm{e}^{-\lambda\, d_{ij}^r}\quad \in \mathbb{R}^{n\times\tilde{n}}, \] which is called \emph{Gibbs kernel} or kernel matrix. 

The following discussion explicitly details the Sinkhorn's algorithm and its properties.

\subsection{Sinkhorn’s Algorithm}
\begin{algorithm}[tbh]
	\scriptsize
		\KwIn{distance  $d_{ij}$ given in \eqref{eq:dist}, probability vectors $p\in\mathbb R_{\geq0}^n$, $\tilde p\in\mathbb R_{\geq0}^{\tilde n}$, regularization parameter $\lambda>0$, $r\ge 1$, threshold $\epsilon$ and  starting value $\tilde\alpha=(\tilde\alpha_1,\dots,\tilde\alpha_{\tilde n})$}
		 Set
		 \begin{equation}\label{eq:MatrixK_sin_div}
			k_{ij}= \exp\big(-\lambda\,d_{ij}^r\big),~\alpha^{(0)} \coloneqq \one_n,\text{ and }\tilde{\alpha}^{(0)}\coloneqq\one_{\tilde n}.
		 \end{equation}
	
		 \While{$\| E^{\Delta} \| > \epsilon $}{
			\If{$\Delta$ is odd}{
				\begin{equation}\hspace*{-7.0cm}
					\begin{split}			
						\alpha_i^{\Delta}\leftarrow&\frac{p_i}{\sum_{j=1}^{\tilde n}k_{ij}\, \tilde\alpha_j^{\Delta-1}},\quad i=1,\dots, n; \\ \tilde{\alpha}_j^{\Delta}\leftarrow &\tilde{\alpha}_j^{\Delta-1}\label{Eq:even_sin_div},\hspace*{1.4cm}\quad j=1,\dots, \tilde n;
					\end{split}
				\end{equation}				
			}
			\Else{
				\begin{equation}\hspace*{-6.5cm}
					\begin{split} 
						\tilde\alpha_j^{\Delta}\leftarrow&\frac{\tilde p_j}{\sum_{i=1}^n\,k_{ij}\,\alpha_i^{\Delta-1}},\quad j=1,\dots, \tilde{n}; \\ \alpha_i^{\Delta}\leftarrow &\alpha_i^{\Delta-1},\label{Eq:odd_sin_div}\hspace*{2.0cm}\quad i=1,\dots, n;
					\end{split}
				\end{equation}}
				increment $\Delta \leftarrow \Delta +1 $}
		 \KwResult{\[s_{r;\lambda}(P,\tilde P) = \frac1\lambda+ \frac1\lambda\sum_{i=1}^n p_i \log\alpha_i^{\Delta^*} +\frac1\lambda\sum_{j=1}^{\tilde n}\tilde p_j\log\tilde\alpha_j^{\Delta^*}-\frac1\lambda\sum_{i=1}^n\sum_{j=1}^{\tilde n} \alpha_i^{\Delta^*}\, \mathrm{e}^{-\lambda\, d_{ij}^r}\, \tilde\alpha_j^{\Delta^*} \]The matrix $\pi^{\Delta^*}=\diag(\alpha^{\Delta^*}) \,  k\,\diag(\tilde{\alpha}^{\Delta^*})$ can also be computed, which is the proximate solution of ERW problem~\eqref{eq:Sinkhorn} \eqref{eq:SD1}--\eqref{eq:SD2}.}
		\caption{Sinkhorn's algorithm \label{alg:Sinkhorn__sin_div}}
	\end{algorithm}
In this section, we illustrate the \emph{iteration process} and stopping criteria of Sinkhorn’s Algorithm~\ref{alg:Sinkhorn__sin_div} to compute ERW distance.

The \emph{iteration counts} of Algorithm~\ref{alg:Sinkhorn__sin_div} are denoted $\Delta \in \mathbb{N}$ and the final iteration count is~$\Delta^*$. 
Algorithm~\ref{alg:Sinkhorn__sin_div} alternately determines $\alpha^\Delta$ and $\tilde \alpha^\Delta$ with
\begin{equation}
	\begin{dcases}
		\tilde{\alpha}^{\Delta} = \tilde{\alpha}^{\Delta-1},& \text{if } \Delta \text{ is odd};\\
		\alpha^{\Delta} = \alpha^{\Delta-1},              & \text{if } \Delta \text{ is even}.\label{eq:OddEvenIteration}
	\end{dcases}
\end{equation}
Sinkhorn's theorem (cf.\ \citet{Sinkhorn1967a, Sinkhorn1967} and Section~\ref{Sec:convergeofSinkhorn} below) for the matrix scaling ensures that iterating~\eqref{eq:S1} converges and the vectors~$\alpha^\Delta$ and~$\tilde\alpha^\Delta$ are unique up to a scalar. From Algorithm~\ref{alg:Sinkhorn__sin_div}, the resultant matrix $\pi^{\Delta^*}=\diag(\alpha^{\Delta^*}) \,  k\,\diag(\tilde{\alpha}^{\Delta^*})$ can be computed, which is the proximate solution of ERW problem~\eqref{eq:Sinkhorn}, \eqref{eq:SD1}--\eqref{eq:SD2}. 

\paragraph{Stopping criteria.}
The iteration process reveals the matrices $\big(\pi^\Delta\big)_{\Delta \in \mathbb{N}}$, which are defined as
\[\pi^\Delta\coloneqq \diag(\alpha^\Delta) \,  k\,\diag(\tilde{\alpha}^\Delta), \quad \Delta \in \mathbb{N}.\]
The norm of residuals $E^\Delta$, which measure the error of the iteration, is
   \begin{equation}\label{eq:Error}
	   \|E^\Delta\| \coloneqq \|\pi^{\Delta}\one_{\tilde n} - p\| + \|(\pi^{\Delta})^{\top}\one_n - \tilde p\|,
   \end{equation} where $\one_n \coloneqq (1,\dots,1)^{\top} \in\mathbb R^n$ and $\one_{\tilde n} \coloneqq (1,\dots,1)^{\top} \in\mathbb R^{\tilde{n}}$.
   If $\Delta$ is odd, $\|(\pi^{\Delta})^{\top}\one_n - \tilde{p}\| = 0$ and if $\Delta$ is even, $\|\pi^{\Delta}\one_{\tilde n} - p\| = 0$.
   The stopping criteria for Algorithm~\ref{alg:Sinkhorn__sin_div}, i.e., $\| E^{\Delta^*} \| \le \epsilon $, implies that
   \begin{equation}
	   \|\pi^{\Delta^*}\one_{\tilde{n}} - p\| + \|(\pi^{\Delta^*})^{\top}\one_n - \tilde{p}\| \le \epsilon.\label{Eq:forThm}
   \end{equation}
Algorithm~\ref{alg:Sinkhorn__sin_div} consolidates the individual steps again. 

\paragraph{Stabilized Sinkhorn’s algorithm.}

The standard Sinkhorn's Algorithm~\ref{alg:Sinkhorn__sin_div} significantly reduces the complexity of the traditional methods.
However, the thirst of larger $\lambda$ among few applications raises the problem of numerical instabilities.
More precisely, for larger $\lambda$, the elementwise exponential matrix $k= \mathrm{e}^{-\lambda\,d^r}$ suffers numerical underflow.
This side effect has increased the need among the OT community to compromise for a slower algorithm, which is known as log-domain stabilized Sinkhorn's algorithm.

The log-domain stabilized Sinkhorn's Algorithm~\ref{alg:Sinkhorn_Log_sin} scales dual variables ($\beta,\gamma$) instead of exponentiated scaling variables ($\alpha,\tilde{\alpha}$), and it utilizes the famous trick among machine learning community called \emph{log-sum-exp trick}.
This log-domain computation and the log-sum-exp trick tackle the numerical underflow. 
\begin{algorithm}[tbh]
	\scriptsize
		\KwIn{distance  $d_{ij}$ given in \eqref{eq:dist}, probability vectors $p\in\mathbb R_{\geq0}^n$, $\tilde p\in\mathbb R_{\geq0}^{\tilde n}$, regularization parameter $\lambda>0$, $r\ge 1$, threshold $\epsilon$ and  starting value $\gamma=(\gamma_1,\dots,\gamma_{\tilde n})$}
		 Set
		 \begin{equation}
			k_{ij}= \exp\big(-\lambda\,d_{ij}^r\big),~\beta^{(0)} \coloneqq \zero_n,\text{ and }{\gamma}^{(0)}\coloneqq\zero_{\tilde n}.
		 \end{equation}
	
		 \While{$\| E^{\Delta} \| > \epsilon $}{
			\If{$\Delta$ is odd}{
				\begin{equation}\hspace*{-7.0cm}
					\begin{split}			
						\hspace*{3.0cm}\beta_i^{\Delta}\leftarrow& \frac1\lambda \Big( \log p_i - \log\big(\sum_{j=1}^{\tilde n}k_{ij}\, \mathrm{e}^{ \lambda \, \gamma_j^{\Delta-1}-\nicefrac12}\big) \Big),\quad i=1,\dots, n; \\ {\gamma}_j^{\Delta}\leftarrow &{\gamma}_j^{\Delta-1},\quad j=1,\dots, \tilde n;\\[-3ex]
					\end{split}
				\end{equation}}				
			\Else{
				\begin{equation}\hspace*{-7.0cm}
					\begin{split} 
						\hspace*{3.0cm}\gamma_j^\Delta\leftarrow&\frac1\lambda\Big( \log \tilde{p}_j - \log\big(\sum_{i=1}^n k_{ij}\, \mathrm e^{ \lambda \, \beta_i^{\Delta-1}-\nicefrac12}\big) \Big),\quad j=1,\dots, \tilde{n}; \\ \beta_i^\Delta\leftarrow &\beta_i^{\Delta-1},\quad i=1,\dots, n;\\[-3ex]
					\end{split}
				\end{equation}}
				increment $\Delta \leftarrow \Delta +1 $}
		 \KwResult{\[s_{r;\lambda}(P,\tilde P) =\sum_{i=1}^np_i\,\beta_i	+ \sum_{j=1}^{\tilde n}\tilde p_j\,\gamma_j -\frac1\lambda\sum_{i=1}^n\sum_{j=1}^{\tilde n} \mathrm e^{-\lambda(d_{ij}^r-\beta_i-\gamma_j)-1}\] The matrix $\pi^{\Delta^*}=\diag({\mathrm{e}^{ \lambda \, \beta^{\Delta^*}-\nicefrac12}}) \,  k\,\diag({\mathrm{e}^{ \lambda \, \gamma^{\Delta^*}-\nicefrac12}})$ can be computed, which is the proximate solution of ERW problem~\eqref{eq:Sinkhorn}, \eqref{eq:SD1}--\eqref{eq:SD2}}
		\caption{Sinkhorn's algorithm (log-domain stabilized) \label{alg:Sinkhorn_Log_sin}}
	\end{algorithm}
	\newpage
Algorithm~\ref{alg:Sinkhorn_Log_sin} encapsulates the individual steps again. 

\subsubsection{Convergence properties of Sinkhorn's iteration}\label{Sec:convergeofSinkhorn}
 The aim of this section is to demonstrate the convergence properties of Sinkhorn's iteration. 
 The following proofs, which summarize the convergence of Sinkhorn's iteration, are applied in many contexts (cf.\ \citet{AltschulerSinkhorn}, \citet{dvurechensky2018computational}, \citet{khalil2018greedy}). 
 We consider the following auxiliary lemmas to substantiate the objective of Algorithm~\ref{alg:Sinkhorn__sin_div} (i.e., the approximation of the Wasserstein distance) from a theoretical standpoint. 
The dual formulation of the ERW problem relates the function $d(\alpha,\tilde\alpha)$ (cf.~\eqref{eq:27}) and
	\begin{equation}
		f(\alpha,\tilde{\alpha})=  \sum_{i=1}^n\sum_{j=1}^{\tilde n} \alpha_i\, \tilde{k}_{ij}\, \tilde\alpha_j -\sum_{i=1}^np_i \log\alpha_i - \sum_{j=1}^{\tilde n}\tilde p_j\log\tilde\alpha_j, \label{Eq:objfunc}
	\end{equation}
	where $\tilde{k}\coloneqq\frac{k}{\|k\|_1}, k=\exp\big(-\lambda\,d^r\big)$, and $p\in\mathbb R_{\geq0}^n$, $\tilde p\in\mathbb R_{\geq0}^{\tilde n} $ are the probability vectors, which satisfy 
	\[ p^{\top}\one_n = \tilde{p}^{\top}\one_{\tilde{n}} = 1.\]
\begin{remark}[Normalization of the kernel matrix $k$]
	In the literature, the approach of normalization of the kernel matrix $k$ is widely used~(cf.\ \citet{AltschulerSinkhorn}, \citet{khalil2018greedy}, \citet{kalantari2008complexity}) for theoretical and numerical analysis. 
	Without loss of generality, we utilize this approach only to substantiate the convergence properties. 
	For numerical experiments, we consider the standard matrix $k$.
\end{remark}
The following Lemma~\ref{lem:Sinkhorniteration1} describes the evolution of the objective function~\eqref{Eq:objfunc} to the target marginals ($p$, $\tilde p$) of Sinkhorn's iteration.
	\begin{lemma}[{cf.\ \citet[Lemma~3.1]{kalantari2008complexity}}]
		The iterates $\alpha^\Delta$ and $\tilde\alpha^\Delta$ of Algorithm~\ref{alg:Sinkhorn__sin_div} satisfy
		\begin{equation}\label{lem:Sinkhorniteration1}
			f(\alpha^\Delta,\tilde\alpha^\Delta)-f(\alpha^{\Delta+1},\tilde{\alpha}^{\Delta+1}) = D_\textit{KL}\big(p\mid\pi^{\Delta}\one_{\tilde n}\big) + D_\textit{KL}\big(\tilde p\mid (\pi^{\Delta})^\top\one_n\big).
		\end{equation}
	\end{lemma}
	\begin{proof}
		First, we assume $\Delta \ge 1$ is even. By equation~\eqref{Eq:objfunc}, it follows that	
		\begin{align}\label{eq:iterationFunctionDifference}
			f(\alpha^{\Delta},\tilde{\alpha}^{\Delta})- f(\alpha^{\Delta+1},\tilde{\alpha}^{\Delta+1})  =  &\sum_{ij}(\alpha_i^\Delta \,\tilde k_{ij}\,\tilde{\alpha}_{j}^{\Delta}-\alpha_{i}^{\Delta+1} \, \tilde k_{ij}\,\tilde\alpha_j^{\Delta+1}) \\
			&+ \sum_{i}p_{i}\big(\log(\alpha_{i}^{\Delta+1})-\log(\alpha_{i}^{\Delta})\big)\\
			& + \sum_j \tilde{p}_j\big(\log(\tilde\alpha_{j}^{\Delta+1})-\log(\tilde{\alpha}_{j}^{\Delta})\big) . 
		\end{align}
		The first component of equation~\eqref{eq:iterationFunctionDifference} turns into  
		\[\sum_{ij}(\alpha_{i}^{\Delta} \,\tilde{k}_{ij}\,\tilde{\alpha}_{j}^{\Delta}-\alpha_{i}^{\Delta+1} \, \tilde{k}_{ij}\,\tilde{\alpha}_{j}^{\Delta+1}) = 0,\] since \[	(\alpha^\Delta)^\top \,\tilde k\,\tilde{\alpha}^{\Delta}= \one_n^{\top}\,\pi^{\Delta}\,\one_{\tilde n}=\tilde{p}^{\top}\,\one_{\tilde n}=1,\] 
		similarly
		\[	(\alpha^{\Delta+1})^\top \,\tilde{k}\,\tilde{\alpha}^{\Delta+1}=\one_{n}^{\top}\,\pi^{\Delta+1}\,\one_{\tilde n}=\one_n^\top\,p=1,\] see \eqref{Eq:even_sin_div}--\eqref{Eq:odd_sin_div}. 
		For this reason, the equation~\eqref{eq:iterationFunctionDifference} becomes  
		\begin{equation}
			\sum_{i}p_{i}\big(\log(\alpha_{i}^{\Delta+1})-\log(\alpha_i^\Delta)\big) + \sum_{j}p_{j}\big(\log(\tilde{\alpha}_{j}^{\Delta+1})-\log(\tilde\alpha_{j}^\Delta)\big).\label{eq:iterationFunctionDifference2}
		\end{equation}By Equation~\eqref{Eq:even_sin_div}, and ~\eqref{Eq:odd_sin_div}, the above Equation~\eqref{eq:iterationFunctionDifference2} becomes
		\[  D_\textit{KL}(p|\pi^{\Delta}\one_{\tilde n}) + D_\textit{KL}(\tilde{p}|(\pi^{\Delta})^{\top}\one_n),\]
		and $D_\textit{KL}(\tilde{p}|(\pi^{\Delta})^{\top}\one_n)=0, \text{ since } \tilde{p}=(\pi^{\Delta})^{\top}\one_n.$ This completes the proof of lemma for $\Delta$ \emph{even} case.
	
		A similar argument applies to the case of \emph{odd}  $\Delta$.
	\end{proof}
In the following Lemma~\ref{lem:Sinkhorniteration2}, we consider the gap between $f(\one_n,\one_{\tilde n})$ and $f(\alpha^{\Delta^*},\tilde{\alpha}^{\Delta^*})$. We know that
	\[f(\one_n,\one_{\tilde{n}}) \coloneqq \ f(\alpha^{(0)},\tilde{\alpha}^{(0)}),\] since $\alpha^{(0)} = \one_n,\text{ and }\tilde{\alpha}^{(0)} = \one_{\tilde n}$, which is a starting value of Algorithm~\ref{alg:Sinkhorn__sin_div}.  
	\begin{lemma}[{cf.\ \citet[Lemma~4.1, Lemma~4.2]{kalantari2008complexity}}]
		\label{lem:Sinkhorniteration2}
		It holds that 
		\begin{equation}
			f(\one_n,\one_{\tilde n})-f(\alpha^{\Delta^*},\tilde{\alpha}^{\Delta^*}) \le \,\log(\nicefrac{\kappa}{\jmath}), 
		\end{equation} where $\kappa$ is the sum of the entries of matrix $\pi^{\Delta^*}$, and $\jmath\coloneqq \min_{ij}\tilde{k}_{ij}.$ 
	\end{lemma}
	\begin{proof}
		Let $(\alpha^{\Delta^*},\tilde{\alpha}^{\Delta^*})$ be the minimizer of the objective function~\eqref{Eq:objfunc}, and we set
		\[\kappa \coloneqq (\alpha^{\Delta^*})^{\top}\,\tilde{k}\,\tilde{\alpha}^{\Delta^*}.\] Equation~\eqref{Eq:objfunc} rewrites as 
		\[f(\alpha^{\Delta^*},\tilde{\alpha}^{\Delta^*}) =  \kappa -\sum_{i=1}^np_i \log\alpha^{\Delta^*}_i - \sum_{j=1}^{\tilde n}\tilde p_j\log\tilde\alpha^{\Delta^*}_j,\] and 
		\[f(\one_n,\one_{\tilde{n}})= \sum_{i=1}^n\sum_{j=1}^{\tilde n} \tilde{k}_{ij}= \kappa.\]
		Now we have
		\begin{equation}
			f(\one_n,\one_{\tilde n})-f(\alpha^{\Delta^*},\tilde{\alpha}^{\Delta^*}) = \sum_{i=1}^np_i \log\alpha^{\Delta^*}_i + \sum_{j=1}^{\tilde n}\tilde p_j\log\tilde\alpha^{\Delta^*}_j.	\label{eq:lem2.10:1}
		\end{equation}
		Without loss of generality we assume that each entry of $\tilde{k}$ is at least $\jmath>0$, then one has
		\begin{equation}\label{eq:lem2.10:3}
			\jmath\Big(\sum_{i=1}^n\alpha^{\Delta^*}_i\Big)\Big(\sum_{j=1}^{\tilde n}\tilde\alpha^{\Delta^*}_j\Big) \le (\alpha^{\Delta^*})^{\top}\,\tilde{k}\,\tilde{\alpha}^{\Delta^*} = \kappa.
		\end{equation}
		Taking the log of both sides of equation~\eqref{eq:lem2.10:3} produces	\begin{equation}
			\log\Big(\sum_{i=1}^n\alpha^{\Delta^*}_i\Big)+\log\Big(\sum_{j=1}^{\tilde n}\tilde\alpha^{\Delta^*}_j\Big) \le \log(\nicefrac{\kappa}{\jmath}). \label{eq:lem2.10:2}
		\end{equation}
		To complete the proof, we consider equations~\eqref{eq:lem2.10:1},~\eqref{eq:lem2.10:2}, and the \label{IterationFinalvalue}log-sum inequality. 
		Now we have 
		\begin{equation}
			\begin{split}
				f(\one_n,\one_{\tilde{n}})-f(\alpha^{\Delta^*},\tilde{\alpha}^{\Delta^*}) = &\sum_{i=1}^np_i \log\alpha^{\Delta^*}_i + \sum_{j=1}^{\tilde n}\tilde p_j\log\tilde\alpha^{\Delta^*}_j \\ \leq \,& \sum_{i=1}^n p_i\,\log\Big(\sum_{l=1}^n\alpha^{\Delta^*}_l\Big)+ \sum_{j=1}^{\tilde{n}} \tilde{p}_j \,\log\Big(\sum_{m=1}^{\tilde n}\tilde\alpha^{\Delta^*}_m\Big) \\ = &\log\Big(\sum_{l=1}^n\alpha^{\Delta^*}_l\Big)+ \log\Big(\sum_{m=1}^{\tilde n}\tilde\alpha^{\Delta^*}_m\Big)  \\ \leq  &\,\log(\nicefrac{\kappa}{\jmath}).
			\end{split} 
		\end{equation}
		This completes the proof of the lemma. 
	\end{proof}

\begin{remark}[Complexity of Sinkhorn’s iteration]\label{rem:sink_iteration_complexity}
	The complexity of Sinkhorn’s iteration is a well studied aspect of regularized Wasserstein problems (cf.\ \citet{AltschulerSinkhorn}, \citet{dvurechensky2018computational}, \citet{khalil2018greedy}).
	Approximately, Sinkhorn's iteration  requires $\mathcal O (\log n +\|d^r\|_\infty\,\lambda)$ arithmetic operations to converge (cf.\ \citet{dvurechensky2018computational}). 
	This means that when $\lambda\to\infty$ corresponding to $\|d^r\|_\infty$ and $n$, number of iteration will be increased.
\end{remark}
\paragraph{Complexity of ERW.} For $n\approx \tilde n$, using Algorithm~\ref{alg:Sinkhorn__sin_div} the ERW problem~\eqref{eq:Sinkhorn}, \eqref{eq:SD1}--\eqref{eq:SD2} can be solved by involving $\mathcal O(n^2 \log n +\|d^r\|_\infty\,\lambda)$ arithmetic operations (note, that $n^2$ operations are needed to perform the matrix-vector multiplication).

\begin{remark}[Entropy bias and Sinkhorn divergence]\label{Sec:sinkhorn_div}

Regardless of the computational advancement of the ERW problem, it is biased.  
That is, 
\begin{equation}
	s_{r;\lambda}(P,P) \neq 0.
\end{equation}
The quantity $s_{r;\lambda}$ is not a distance, more specifically, it violates the axiom of definiteness of the distance function.
To overcome this difficulty, \citet{ramdas2017wasserstein} introduce the \emph{Sinkhorn divergence}  as
\begin{equation}
	{sd}_{r;\lambda}(P,\tilde P) \coloneqq  s_{r;\lambda}(P,\tilde P) - \frac{1}{2}s_{r;\lambda}(P,P) -\frac{1}{2}s_{r;\lambda}(\tilde{P},\tilde{P}),
\end{equation}
which is a natural normalization (or debias) of the quantity.
The key properties of Sinkhorn divergence include
\begin{enumerate}[nolistsep]
	\item non-negativity,
	\item $\lim_{\lambda\to\infty}{sd}_{r;\lambda}(P,\tilde P)= w_r(P,\tilde P)$  and
	\item ${sd}_{r;\lambda}(P,P)=0$ for all $\lambda>0$.
\end{enumerate}
\end{remark}

\section{Nonequispaced Fast Fourier Transform~(NFFT)}
\label{sec:Fast_summation}
Generally, for a faster computation of the matrix-vector multiplication with the distance matrix~$d^r$, equispaced convolution is used.
The most common algorithm used to compute equispaced convolution is the standard FFT algorithm.
Our research, in contrast, promotes the \emph{non}equispaced convolution, which is approximated by the \emph{non}equispaced fast Fourier transform (NFFT) to accelerate the computation of the Sinkhorn's Algorithm~\ref{alg:Sinkhorn__sin_div}.
More precisely, the matrix-vector multiplications of Sinkhorn’s iteration~\eqref{eq:S1}, which is the main computational bottleneck, are tackled by \emph{fast summation} based on NFFT in $\mathcal O\big(n \log n\big)$ arithmetic operations. 
Moreover, this fast summation technique has better stability and is accomplished with machine precision.

\subsection{NFFT-based fast summation}
This subsection succinctly describes the fast summation based on NFFT. 

The fast summation method based on NFFT takes advantage of the special structure of the \emph{Euclidean} distance matrix. 
The distance matrix $d \in \mathbb{R}^{n\times\tilde{n}}$ in~\eqref{eq:Da} has entries
\begin{equation}\label{eq:dist}
	d(x_i,\tilde x_j)\coloneqq \| x_i- \tilde{x}_j\|,
\end{equation}
which is the distances of all combinations of states, and we recall that $\|\cdot\|$ denotes the Euclidean norm or $2$-norm.

\paragraph{Approximation of matrix-vector multiplication of Sinkhorn's iteration.}
The fast summation technique based on NFFT takes advantage of the particular form of the sums
\begin{equation}\label{eq:matexp}
	t_i \colonapprox \big(k\,\tilde{\alpha}\big)_{i} =\sum_{j=1}^{\tilde n} \tilde \alpha_j\, \mathrm{e}^{-\lambda \| x_i -\tilde x_j\|^r},\quad i=1,\dots, n,
\end{equation}
as well as of the sums of the ‘transposed’,
\begin{equation}\label{eq:matexpt}
	\tilde t_j\colonapprox \big(k\,{\alpha}\big)_{j}=\sum_{i=1}^{n} \alpha_i\, \mathrm{e}^{-\lambda \| x_i -\tilde x_j\|^r},\quad j=1,\dots,\tilde n,
\end{equation}
since these summations are the bottleneck of Sinkhorn's iteration.
\paragraph{An overview of NFFT.}
For \emph{equispaced} points $x_i$ and $\tilde x_j$, the summation of \eqref{eq:matexp} and \eqref{eq:matexpt} correspond to the  multiplication of a Toeplitz matrix with a vector, respectively.
In this case, we immediately obtain a fast algorithm based on embedding the matrix into a circulant matrix and then diagonalize the matrix by the Fourier matrix, see \citet[Theorem~3.31]{PlPoStTa18}, such that we end up with ${\cal O}(n\log n)$ operations using the FFT, see \citet[pp.\ 141--142]{PlPoStTa18}.
A fast algorithm with \emph{arbitrary} points follows the similar ideas, but based on the NFFT, see \citet[Chapter~7]{PlPoStTa18} and the related software in \citet{nfft3} for details. 

\subsubsection{The ansatz of NFFT based fast summation}
The core idea of \emph{fast summation} based on NFFT is to accurately approximate the radial kernel function 
\begin{equation}\label{eq:kernel_fun} \mathcal K(y) \coloneqq  {\mathrm e}^{-\lambda \| y\|^r}. \end{equation}
In general, this approximation of the kernel function $\mathcal K(y)$ accommodates, when the entries of the matrix $k$ are in the form  \begin{equation}\label{eq:kernel_fun_entries}k_{ij}= \mathcal K(x_i-\tilde x_j).\end{equation}

In terms of fast summation, the goal of the NFFT is to accurately approximate $\mathcal K(y)$ by a $h$-periodic trigonometric polynomial $\mathcal K_{RK}(y)$,
	\begin{equation}\label{eq:ansatz_NFFT}
		\mathcal K(y) \approx \mathcal K_{RK}(y) \coloneqq \sum_{\mathrm k \in \mathcal{I}_N} b_\mathrm k\, \mathrm e^{2  \pi \mathrm i \mathrm k y / h}, \quad \mathcal{I}_N \coloneqq \Big\{ -\frac{N}{2},-\frac{N}{2}+1,\dots,-1,0,\dots,\frac{N}{2}-1\Big\}^d,
	\end{equation}
	with appropriate Fourier coefficients $b_{\mathrm k}\in\mathbb C$ and bandwidth $N \in 2\mathbb{N}$.
For example, when  considering the Gaussian kernel function~($r=2$), we have
\begin{equation}\label{eq:ex_gauss_kernel_1} \big(k\,{\alpha}\big)_{j} \ =\sum_{i=1}^{n} \alpha_i\, \mathrm{e}^{-\lambda \| x_i -\tilde x_j\|^2},\quad  j=1,\dots,\tilde n.
\end{equation} 
Now we rewrite Equation~\eqref{eq:ex_gauss_kernel_1} by involving the kernel function $\mathcal K(y) = \mathrm  {\mathrm e}^{-\lambda \| y\|^2}$ as 
\begin{equation}\label{eq:ex_gauss_kernel_2} \big(k\,{\alpha}\big)_{j}   \coloneqq \sum_{i=1}^{n} \alpha_i\,  \mathcal K(x_i -\tilde x_j),\quad  j=1,\dots,\tilde n.
\end{equation}
For the efficient computation of \eqref{eq:ex_gauss_kernel_2}, the fast summation technique based on NFFT approximates~$\mathcal K$ by the trigonometric polynomial $\mathcal K_{RK}$.

\begin{figure}[!ht]
	\centering
		\begin{tikzpicture}[scale=0.8]
			\begin{axis}[
			width=11cm,
			height=5cm,
			axis x line=none,
			axis y line=none,
			xmin=-2,xmax=4.5,
			ymin=-0.5,ymax=1.2]
			
			\def\x{-0.2};
			
			\addplot[color=black,line width=0.8pt] coordinates {(-1.5,\x) (4,\x)};
			\addplot[smooth,line width=0.8pt,restrict x to domain=-1.25:-0.75,densely dotted] table[x=x,y expr=\thisrow{y}] {kern.txt};
			\addplot[smooth,line width=0.8pt,restrict x to domain=-0.75:0.75] table[x=x,y expr=\thisrow{y}] {kern.txt};
			\addplot[smooth,line width=0.8pt,restrict x to domain=0.75:1.75,densely dotted] table[x=x,y expr=\thisrow{y}] {kern.txt};
			\addplot[smooth,line width=0.8pt,restrict x to domain=1.75:3.25] table[x=x,y expr=\thisrow{y}] {kern.txt};
			\addplot[smooth,line width=0.8pt,restrict x to domain=3.25:3.75,densely dotted] table[x=x,y expr=\thisrow{y}] {kern.txt};
			
			\addplot[line width=0.8pt,mark=|] coordinates {(-1.25,\x)} node[below] {$-\frac h2$};
			\addplot[line width=0.8pt,mark=|] coordinates {(-0.75,\x)} node[below] {\small$-L$};
			\addplot[line width=0.8pt,mark=|] coordinates {(0.75,\x)} node[below] {\small$L$};
			\addplot[line width=0.8pt,mark=|] coordinates {(1.25,\x)} node[below] {$\frac h2$};
			\addplot[line width=0.8pt,mark=|] coordinates {(1.75,\x)} node[below] {\small$h\!-\!L$};
			\addplot[line width=0.8pt,mark=|] coordinates {(3.75,\x)} node[below] {$\frac{3h}2$};
			
			\addplot[line width=0.8pt,mark=*] coordinates {(-0.75,0.2821)};
			\addplot[line width=0.8pt,mark=*] coordinates {(0.75,0.2821)};
			
			\addplot[line width=0.8pt,mark=*] coordinates {(1.25,0.04)};
			
			\addplot[line width=0.8pt,mark=*] coordinates {(1.75,0.2821)};
			\addplot[line width=0.8pt,mark=*] coordinates {(3.25,0.2821)};
			
			\node (A) at (axis cs:1.15,0.5){};
			\node (B) at (axis cs:0.75,0.3){};
			\draw[->] (A)--(B);
			
			\node (A) at (axis cs:1.25,0.45){};
			\node (B) at (axis cs:1.25,0.06){};
			\draw[->] (A)--(B);
			
			\addplot[color=black,mark=text,/pgf/text mark={\text{smooth}}] coordinates{(1.4,0.53)};
			
			\addplot[color=black,mark=text,/pgf/text mark={\small$\mathcal K(\cdot)$}] coordinates{(0,1.1)};
			\addplot[color=black,mark=text,/pgf/text mark={\footnotesize$\mathcal K_{\rm B}(\cdot)$}] coordinates{(0.75,0.0)};
			
			\end{axis}
		\end{tikzpicture}	
		\begin{tikzpicture}[scale=0.8]
			\begin{axis}[
			width=10cm,
			height=10cm,
			axis x line=none,
			axis y line=none,
			xmin=-2.5,xmax=4,
			ymin=-0.5,ymax=6,
			set layers,
			mark layer=axis tick labels
			]
			\addplot[color=black,line width=0.8pt] coordinates {(-1.75,0) (4,0)};
			\addplot[color=black,line width=0.8pt] coordinates {(-1.75,0) (-1.75,3)};
			\addplot[line width=0.8pt,mark=|] coordinates {(-0.75,0)} node[below] {\small$-L$};
			\addplot[line width=0.8pt,mark=|] coordinates {(0.75,0)} node[below] {\small$L$};
			\addplot[line width=0.8pt,mark=|] coordinates {(-1.25,0)} node[below] {\tiny $-\frac h2$};
			\addplot[line width=0.8pt,mark=|] coordinates {(1.25,0)} node[below] {\tiny $\frac h2$};
			\addplot[line width=0.8pt,mark=|] coordinates {(3.75,0)} node[below] {\tiny $\frac{3h}2$};
			\addplot[line width=0.8pt,mark=|] coordinates {(1.75,0)} node[below] {\small$h\!-\!L$};
			\addplot[line width=0.8pt,mark=-] coordinates {(-1.75,-0.75+1.5)} node[left] {\small$-L$};
			\addplot[line width=0.8pt,mark=-] coordinates {(-1.75,0.75+1.5)} node[left] {\small$L$};
			\addplot[line width=0.8pt,mark=-] coordinates {(-1.75,-1.25+1.5)} node[left] {$-\frac h2$};
			\addplot[line width=0.8pt,mark=-] coordinates {(-1.75,1.25+1.5)} node[left] {$\frac h2$};
			\addplot[color=black,mark=square*,/tikz/mark size=45.9,line width=1pt,mark options={fill=black!25}] coordinates {(0,1.5)};
			\addplot[color=black,mark=*,/tikz/mark size=45.75,line width=1pt,mark options={fill=black!15}] coordinates {(0,1.5)};
			\addplot[color=black,mark=*,/tikz/mark size=27.5,line width=1pt,mark options={fill=white}] coordinates {(0,1.5)};
			\addplot[color=black,mark=square*,/tikz/mark size=45.9,line width=1pt,mark options={fill=black!25}] coordinates {(2.5,1.5)};
			\addplot[color=black,mark=*,/tikz/mark size=45.75,line width=1pt,mark options={fill=black!15}] coordinates {(2.5,1.5)};
			\addplot[color=black,mark=*,/tikz/mark size=27.5,line width=1pt,mark options={fill=white}] coordinates {(2.5,1.5)};
			\addplot[color=black,line width=0.8pt,dashed] coordinates {(-1.75,0.25) (4,0.25)}; 
			\addplot[color=black,line width=0.8pt,dashed] coordinates {(-1.75,0.75) (4,0.75)}; 
			\addplot[color=black,line width=0.8pt,dashed] coordinates {(-1.75,2.25) (4,2.25)}; 
			\addplot[color=black,line width=0.8pt,dashed] coordinates {(-1.75,2.75) (4,2.75)}; 
			\addplot[color=black,line width=0.8pt,dashed] coordinates {(-1.25,0) (-1.25,3)}; 
			\addplot[color=black,line width=0.8pt,dashed] coordinates {(-0.75,0) (-0.75,3)}; 
			\addplot[color=black,line width=0.8pt,dashed] coordinates {(0.75,0) (0.75,3)}; 
			\addplot[color=black,line width=0.8pt,dashed] coordinates {(1.25,0) (1.25,3)}; 
			\addplot[color=black,line width=0.8pt,dashed] coordinates {(1.75,0) (1.75,3)}; 
			\addplot[color=black,line width=0.8pt,dashed] coordinates {(3.25,0) (3.25,3)}; 
			\addplot[color=black,line width=0.8pt,dashed] coordinates {(3.75,0) (3.75,3)}; 
			\addplot[color=black,mark=text,/pgf/text mark={\footnotesize $\mathcal K(\| y\|)$}] coordinates{(0,1.5)};
			\addplot[color=black,mark=text,/pgf/text mark={\footnotesize $\mathcal K_{\rm B}(\|y\|)$}] coordinates{(0,2.45)};
		\end{axis}
		\end{tikzpicture}
	\caption{The regularized periodic function --- dimension one (left) and the dimensions two (right), see \citet[Figure 3.14]{nestlerdiss}.}\label{Fig:per}
\end{figure}
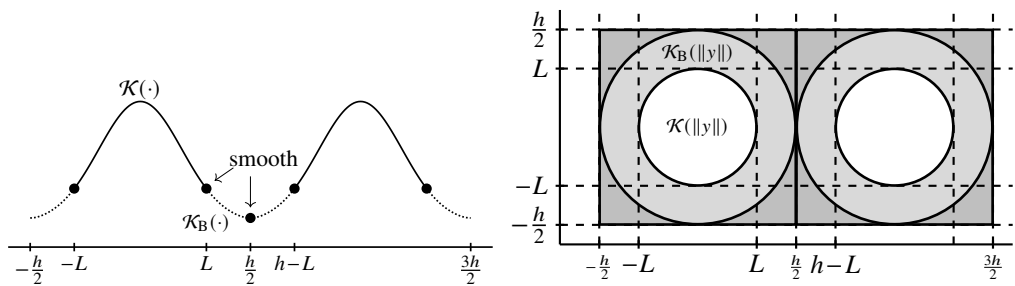

From Equation~\eqref{eq:ansatz_NFFT}, we notice that $\mathcal K_{RK}(y)$ are $h$-periodic functions, although the kernel $\mathcal K(y)$ is not $h$-periodic.
Therefore, we regularize ${\mathcal K}(y)$ to obtain a $h$-periodic smooth kernel function $\tilde{\mathcal K}(y)$, which is $p-1$ times continuously differentiable in the periodic setting, where $p \in \mathbb{N}$ is the degree of smoothness, and the Fourier coefficients decay quickly.
\paragraph{Regularization of ${\mathcal K}(y)$.}
Assume that we have $\|x_j\|\leq \frac L2$, i.e.,  $\|x_i-\tilde x_j\|\leq L$, for some $L>0$.
We define the  multivariate, $h$-periodic regularized kernel function $\tilde{\mathcal K}\colon [-\frac{h}{2},\frac{h}{2}]^d\to \mathbb R$ with $h\ge 2L$ by
\[
	\tilde{\mathcal K}(y)\coloneqq
		\begin{cases}
		\mathcal K(\|y\|) &\text{if } \|y\|\leq L, \\
		\mathcal K_\mathrm{B}(\|y\|) &\text{if } L<\|y\|\leq\frac h2, \\
		\mathcal K_\mathrm{B}(\frac h2) &\text{if } y\in[-\frac h2,\frac h2]^d\text{ and } \|y\|>\frac h2,
		\end{cases}
\]
where $K_\mathrm{B}$ is an appropriately chosen univariate polynomial, which is constructed using two-point Taylor interpolation, see Figure~\ref{Fig:per}. For a detailed interpretation of this approach, we refer to \citet[Chapter~7.5]{PlPoStTa18}. 
\paragraph{Approximation of smooth periodic function $\tilde{\cal  K}$.}
In the univariate case, we are now able to approximate the smooth periodic function $\tilde{\cal  K}$ by a Fourier series to obtain
\begin{equation}\tilde t_j \coloneqq \sum_{i=1}^{n} \alpha_i\, \tilde{\mathcal{K}}(x_i -\tilde x_j) \approx \sum_{i=1}^{n} \alpha_i\, {\mathcal{K}_{RK}}(x_i -\tilde x_j) ,\quad j=1,\dots,\tilde n. \end{equation}	
By using~\eqref{eq:ansatz_NFFT} and interchanging the order of summation, as well as utilizing the outstanding property 
\[\mathrm{e}^{2\pi \mathrm{i} (x_i-\tilde x_j)/h}=\mathrm{e}^{2\pi \mathrm i x_i/h}\, \mathrm{e}^{-2\pi \mathrm{i} \tilde x_j/h},\] we obtain
\begin{equation}\label{eq:fastsum}
	t_i\approx 
	\sum_{\mathrm k \in \mathcal{I}_N} b_\mathrm k \left(
	\sum_{j=1}^{\tilde n} \tilde \alpha_j \mathrm e^{-2\pi \mathrm i \mathrm k \tilde x_j/h}
	\right)
	\mathrm{e}^{2\pi \mathrm{i} \mathrm k x_i/h}, \quad i=1,\dots, n.
\end{equation}
We compute the inner sum for each $\mathrm k \in \mathcal{I}_N$ using the NFFT in $\mathcal O(N\log N +\tilde n)$ arithmetical operations and the outer sum with $\mathcal O(N\log N + n)$.

This simple idea works very well, if the function $\mathcal K$ is smooth and can be approximated by a short Fourier series $\tilde{\mathcal K}$,i.e., by a small polynomial of degree $N$.
This is especially true for the case $r=2$, where the method is also known as \emph{fast Gaussian transform}. 
%
\begin{remark}
	We note as well that large values of~$\lambda$ corresponds to a localization of the support points of the measure. In this setting, the set of support points can be ordered in reduced operations (as mentioned below), so that matrix-vector operations are eligible in the same time as our implementation.
	For this reason, our algorithm is primarily adapted for small values of~$\lambda$. Nevertheless, it renders stable computation for sufficiently large values of~$\lambda$.
\end{remark}
\begin{remark}[Arithmetic complexity]\label{rem:complexity}
	For $\lambda>0$, the kernel approximation~\eqref{eq:ansatz_NFFT} is independent of $n$ ($\tilde{n},$ resp.) data points, therefore we can appropriately fix the  polynomial degree $N$. Thus, the approximation   ends up with $\mathcal O (n + \tilde{n})$ arithmetic operations. 
	Furthermore, for $r=1$, we need additional near-field regularisation at the point 
	$y=0$.
	In this case, we end up with an arithmetic complexity of $\mathcal O (n \log n+ \tilde{n}\log \tilde{n})$, see \citet[Chapter~7.5]{PlPoStTa18} for a detailed interpretation.
\end{remark}
\begin{figure}[!htb]
	\begin{center}
		\begin{tikzpicture}[scale=0.5]
			\begin{axis}[
				xlabel={$y$},
				xmin=-3, xmax=3,
				ymin=-0.1, ymax=1.1,
				legend pos=north east,
				]
				\addplot[color=violet, samples=500]{exp(-x*x)};	
				\legend{exp(-$y^2$),exp(-$|y|$)
				}
				\addplot[color=green, samples=500]{exp(-abs(x))};

				\end{axis}
				\end{tikzpicture}
				
		\caption{Wasserstein Distance for $r=1$ (green) and $r=2$ (violet).}\label{Fig:Dist}
	\end{center}
\end{figure}
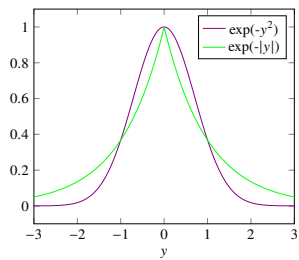

\subsection{NFFT boost for the Sinkhorn's algorithm\label{Sec:FFT_boost}}
This section presents the NFFT-accelerated standard and log-domain Sinkhorn's algorithms.
The NFFT-accelerated Sinkhorn's algorithms propose a novel method using non-equispaced convolution to approximate the Wasserstein distance. 
The algorithms 
below describe the operations of our proposed method schematically. 

The NFFT-accelerated ERW distance (lower bound) including the entropy is denoted by
		\[
			s_{r;\lambda;\text{NFFT}}(P,\tilde P),
		\]
which can be computed using Algorithm~\ref{alg:SinkhornFFT} and \ref{alg:logSinkhornFFT}. This quantity is an approximation of  $s_{r;\lambda}(P,\tilde P)$, i.e., $s_{r;\lambda}(P,\tilde P) \colonapprox s_{r;\lambda;\text{NFFT}}(P,\tilde P)$. Furthermore, the ERW distance (upper bound) is computed by 
		\[
			\tilde s_{r;\lambda;\text{NFFT}}(P,\tilde P) = s_{r;\lambda;\text{NFFT}}(P,\tilde P)+ \frac{H(P)+H(\tilde P)}\lambda. 
		\]
The NFFT-accelerated Sinkhorn divergence is computed by  \[{sd}_{r;\lambda;\text{NFFT}}(P,\tilde P)\coloneqq  s_{r;\lambda}(P,\tilde P) - \frac{1}{2}s_{r;\lambda;\text{NFFT}}(P,P) -\frac{1}{2}s_{r;\lambda;\text{NFFT}}(\tilde{P},\tilde{P}).\]
\paragraph{Arithmetic complexity.} For simplicity, we assume that $n=\tilde n$. As mentioned eariler, the evaluation of sums in Algorithm~\ref{alg:SinkhornFFT} and~\ref{alg:logSinkhornFFT} take only $\mathcal O (n)$ arithmetic operations for $r=2$ and $\mathcal O (n \log n)$ for $r=1$. From Remark~\ref{rem:sink_iteration_complexity}, we know that Sinkhorn's iteration process requires $\mathcal O (\log n +\|d^r\|_\infty\,\lambda)$. Therefore, for $r=2$, our proposed algorithms require only
\[ \mathcal O (n \, \log n +\|d^r\|_\infty\,\lambda), \] and for $r=1$ \[ \mathcal O (n \,(\log n)^2 +\|d^r\|_\infty\,\lambda). \]
\begin{algorithm}[!tbh]
	\scriptsize
		\KwIn{support nodes $x_i$ and $\tilde x_j$, probability vectors $p\in\mathbb R_{\geq0}^n$, $\tilde p\in\mathbb R_{\geq0}^{\tilde n}$, regularization parameter $\lambda>0$, threshold $\varepsilon $ and  starting value $\tilde\alpha=(\tilde\alpha_1,\dots,\tilde\alpha_{\tilde n})$}
		\While{$\| E^{\Delta} \| > \varepsilon $}
		{	Set
		\begin{equation}\label{eq:starting_sinkhorn_NFFT}
		   \alpha^{(0)} \coloneqq \one_n,\text{ and }\tilde{\alpha}^{(0)}\coloneqq\one_{\tilde n}.
		\end{equation}
			
			\If{$\Delta$ is odd}{compute
			\[	t_i^{\Delta-1}\leftarrow \sum_{j=1}^{\tilde n} \tilde\alpha_j^{\Delta-1}\, \mathrm{e}^{-\lambda \|x_i-\tilde x_j\|^r},\quad i=1,\dots, n,\] by employing the fast summation~\eqref{eq:fastsum} and set
			\begin{equation}\hspace*{-7.0cm}\begin{split} \alpha_i^{\Delta}\leftarrow&\frac{p_i}{t_i^{\Delta-1}},\quad i=1,\dots, n; \\ \tilde{\alpha}_j^{\Delta}\leftarrow &\tilde{\alpha}_j^{\Delta-1},\label{Eq:even_F_sink}\hspace*{-0.1cm}\quad j=1,\dots, \tilde{n}.\\[-2ex]\end{split}\end{equation}				
			}
			\Else{compute 
				\[\tilde t_j^{\Delta-1}\leftarrow \sum_{i=1}^n  \mathrm{e}^{-\lambda \|x_i-\tilde x_j\|^r}\,\alpha_i^{\Delta-1},\quad j=1,\dots,\tilde n,\] by employing the fast summation  \eqref{eq:fastsum} and set
				\begin{equation}\hspace*{-7.0cm}\begin{split} \tilde\alpha_j^{\Delta}\leftarrow&\frac{\tilde p_j}{\tilde t_j^{\Delta-1}}, \quad j=1,\dots, \tilde{n}; \\ \alpha_i^{\Delta}\leftarrow &\alpha_i^{\Delta-1},\label{Eq:odd_F_sink}\hspace*{-0.1cm}\quad i=1,\dots, n.\\[-2ex]\end{split}\end{equation}}
				increment $\Delta \leftarrow \Delta +1$
			}
		\KwResult{The ERW distance (cf.~\eqref{eq:14}) approximating the Wasserstein distance $W_r(P,\tilde P)$ is
		\[	s_{r;\lambda;\text{NFFT}}(P,\tilde P)\coloneq \frac1\lambda+\frac1\lambda\sum_{i=1}^np_i \log\alpha_i^{\Delta^*} +\frac1\lambda\sum_{j=1}^{\tilde n}\tilde p_j\log\tilde\alpha_j^{\Delta^*}-\frac1\lambda\sum_{j=1}^{\tilde n} \tilde t_j^{\Delta^*}\, \tilde\alpha_j^{\Delta^*}.\] }
		\caption{NFFT-accelerated Sinkhorn's algorithm\label{alg:SinkhornFFT}}
	\end{algorithm}
	\begin{algorithm}[!tbh]
		\scriptsize
			\KwIn{support nodes $x_i$ and $\tilde x_j$, probability vectors $p\in\mathbb R_{\geq0}^n$, $\tilde p\in\mathbb R_{\geq0}^{\tilde n}$, regularization parameter $\lambda>0$, threshold $\varepsilon $ and  starting value $\gamma=(\gamma_1,\dots,\gamma_{\tilde n})$}
			\While{$\| E^{\Delta} \| > \varepsilon $}
			{	Set
			\begin{equation}\label{eq:starting_Log_sinkhorn_NFFT}
			   \beta^{(0)} \coloneqq \zero_n,\text{ and }{\gamma}^{(0)}\coloneqq\zero_{\tilde n}.
			\end{equation}
				
				\If{$\Delta$ is odd}{compute
				\[	t_i^{\Delta-1}\leftarrow \sum_{j=1}^{\tilde n} \mathrm{e}^{ \lambda \, \gamma_j^{\Delta-1}-\nicefrac12}\, \mathrm{e}^{-\lambda \|x_i-\tilde x_j\|^r},\quad i=1,\dots, n,\] by employing the fast summation~\eqref{eq:fastsum} and set
				\begin{equation}\hspace*{-5.5cm}\begin{split} \beta_i^{\Delta}\leftarrow&\frac{1}{\lambda} \Big( \log p_i - \log t_i^{\Delta-1} \Big) ,\quad i=1,\dots, n; \\ \gamma_j^{\Delta}\leftarrow &\gamma_j^{\Delta-1},\label{Eq:even_F_sink_log}\hspace*{-0.1cm}\quad j=1,\dots, \tilde{n}.\\[-2ex]\end{split}\end{equation}				
				}
				\Else{compute 
					\[\tilde t_j^{\Delta-1}\leftarrow \sum_{i=1}^n  \mathrm{e}^{-\lambda \|x_i-\tilde x_j\|^r}\,\mathrm{e}^{ \lambda \, \beta_i^{\Delta-1}-\nicefrac12},\quad j=1,\dots,\tilde n,\] by employing the fast summation  \eqref{eq:fastsum} and set
					\begin{equation}\hspace*{-5.5cm}\begin{split} \gamma_j^{\Delta}\leftarrow&\frac{1}{\lambda}\Big( \log \tilde{p}_j - \log\tilde t_j^{\Delta-1} \Big), \quad j=1,\dots, \tilde{n}; \\ \beta_i^{\Delta}\leftarrow &\beta_i^{\Delta-1},\label{Eq:odd_F_sink_log}\hspace*{-0.1cm}\quad i=1,\dots, n.\\[-2ex]\end{split}\end{equation}}
					increment $\Delta \leftarrow \Delta +1$
				}
			\KwResult{The ERW distance (cf.~\eqref{eq:14}) approximating the Wasserstein distance $W_r(P,\tilde P)$ is
			\[	s_{r;\lambda}(P,\tilde P)\approx \sum_{i=1}^np_i\,\beta_i	+ \sum_{j=1}^{\tilde n}\tilde p_j\,\gamma_j -\frac1\lambda\sum_{j=1}^{\tilde n} \tilde t_j^{\Delta^*}\,\mathrm{e}^{ \lambda \, \gamma_j^{\Delta^*}-\nicefrac12}.\] }
			\caption{NFFT-accelerated Sinkhorn's algorithm (log-domain)\label{alg:logSinkhornFFT}}
		\end{algorithm}
\begin{remark}[Optimal transition matrix $\pi^{\Delta^*}$]
	The NFFT-accelerated Sinkhorn's Algorithm~\ref{alg:SinkhornFFT} and \ref{alg:logSinkhornFFT} bypass the allocations of the matrices $d$, $k$ and $\pi^{\Delta^*}$ and returns the objective of the Sinkhorn's algorithm, i.e., the ERW distance $s_{r;\lambda}(P,\tilde P)$ of the measures~$P$ and~$\tilde P$.
	Our proposed algorithms provide the optimal exponentiated  dual variables $(\alpha,\tilde{\alpha})$~(Algorithm~\ref{alg:SinkhornFFT}) and optimal dual variables $(\beta^{\Delta^*},\gamma^{\Delta^*})$~(Algorithm~\ref{alg:logSinkhornFFT}). Hence, the optimal transition matrix $\pi^{\Delta^*}$ can still be computed with~\eqref{eq:pi*}. However, this -- as mentioned~-- requires $\mathcal O(n\,\cdot\tilde n)$ operations, which would increase the performance time and thus is avoided.
\end{remark}

The NFFT fast summation technique splendidly adapts to the Sinkhorn’s algorithms. As mentioned eariler, this technique guarantees fast and memory efficient computation with  machine accuracy.


\section{Numerical Experiments}
\label{Sec:Numerical_exposition}
This section demonstrates the performance and accuracy of NFFT-accelerated Sinkhorn's Algorithm~\ref{alg:SinkhornFFT} using synthetic as well as real data sets. 
All runtime measurements were performed on a standard desktop computer with Intel(R) Core(TM) i7-7700 CPU and 15.0 GB of RAM.
The source code of our implementation of standard Sinkhorn's Algorithm~\ref{alg:Sinkhorn__sin_div}, log-domain stabilized Sinkhorn's Algorithm~\ref{alg:Sinkhorn_Log_sin}, NFFT accelarated Sinkhorn's Algorithm~\ref{alg:SinkhornFFT}, NFFT accelarated log-domain Sinkhorn's Algorithm~\ref{alg:logSinkhornFFT} and linear programming solver to compute $w_r(P,\tilde P)$, which can be used to reproduce the following results, are available in online.\footnote{Cf.\ \href{https://github.com/rajmadan96/NFFT-Sinkhorn-Wasserstein_distance/}{https://github.com/rajmadan96/NFFT-Sinkhorn-Wasserstein$\_$distance/}} The implementation of our proposed algorithms are based on the freely available repository ‘NFFT3.jl’.\footnote{Cf.\ \href{https://github.com/NFFT/NFFT3.jl}{https://github.com/NFFT/NFFT3.jl}}
%
\subsection{Synthetic data} 
\label{Sec:Synthetic data}
We test in one-dimension (Section~\ref{sec:17} below), and  for two-dimensional data (Section~\ref{sec:18}) to  demonstrate the performance of our proposed algorithm which still delivers results,  which are out of reach for tradtional implementations.
\subsubsection{NFFT-accelerated Sinkhorn's algorithm in one dimension}\label{sec:17}
Consider a measure~$P$ on~$\mathbb R$ with quantiles $s_i$, i.e., 
\[	P\big((-\infty,s_i]\big)= \frac i{\tilde n+1}, \qquad i=1,\dots,\tilde n,\]
and corresponding weights 
\[	p_i\coloneqq P\left(\left(\frac{s_{i-1}+s_i}2,\frac{s_i+s_{i+1}}2\right]\right),\qquad i=1,\dots, \tilde n.\]
The measure
\[\tilde P_{\tilde n}\coloneqq\sum_{i=1}^{\tilde n} p_i\, \delta_{s_i} \]
is the best discrete approximation of~$P$ in Wasserstein distance (cf.\ \citet{GrafLuschgy}).

To demonstrate the performance of Algorithm~\ref{alg:SinkhornFFT}, we consider independent and identically distributed observations $X_i\in\mathbb R$, $i=1,\dots,n$, from the measure~$P$, and the corresponding empirical measure
\[ \hat P_n\coloneqq \frac1n\sum_{i=1}^n \delta_{X_i}. \]

Table~\ref{tab:Sinkhorn} compares the  computation time of Sinkhorn's Algorithm~\ref{alg:Sinkhorn__sin_div}, and the NFFT-accelerated Sinkhorn's Algorithm~\ref{alg:SinkhornFFT}.

\begin{table}[!ht]
	\begin{centering}
		\begin{tabular}{lcccccc}
			\toprule 
			$n=\tilde n$:  & \num{1000} & \num{10000}&\num{10000} & \num{100000} &\num{1000000}&\num{10000000}\tabularnewline
			\midrule 
			Sinkhorn's Algorithm~\ref{alg:Sinkhorn__sin_div} 
						 & \num{.49}\,s & \num{23.72}\,s & \num{41.16}\,s&\multicolumn{3}{c}{\emph{out of memory or > 1\,hour}}\tabularnewline[3\doublerulesep]
			NFFT-accelerated Sinkhorn~\ref{alg:SinkhornFFT} 
						 & 0.28\,s &  0.39\,s &2.08\,s & 2.31\,s & 9.38\,s & 62.4\,s\tabularnewline[3\doublerulesep]
						
			\bottomrule
		\end{tabular}
	\par\end{centering}
	\caption{\label{tab:Sinkhorn}Dimension~1: Comparison of computation times for $r=2$ and $\lambda=20$}
\end{table}

The table demonstrates that the NFFT-accelerated Sinkhorn's Algorithm~\ref{alg:SinkhornFFT} easily delivers results for problem sizes, which are out of reach for the traditional Sinkhorn's Algorithm~\ref{alg:Sinkhorn__sin_div}.

\subsubsection{NFFT-accelerated Sinkhorn's in two dimension}\label{sec:18}

We demonstrate next performance of the NFFT-accelerated Sinkhorn's Algorithm~\ref{alg:SinkhornFFT} by approximating the Wasserstein distance for empirical measures
\[	P=\frac1n\sum_{i=1}^n\delta_{(U^1_i,U^2_i)}\
	\text{ and }\
	\tilde P=\frac1{\tilde n}\sum_{j=1}^{\tilde n}\delta_{(\tilde U_j^1,\tilde U_j^2)}\]
on $\mathbb R\times\mathbb R$, where $U^1_i$, $U^2_i$, $i=1,\dots,n$, and $\tilde U_j^1$, $\tilde U_j^2$, $j=1,\dots,\tilde n$, are independent samples from the uniform distribution.

Table~\ref{tab:Sinkhorn2} displays execution times for the uniform distribution on $[0,1]\times [0,1]$.
While computation time and memory allocations are already critical for $n$, $\tilde n\approx \num{100000}$, the NFFT-accelerated Sinkhorn's algorithm still performs in reasonable time.

\begin{table}[htb]
	\begin{centering}
		\begin{tabular}{lccccc}
			\toprule 
			$n=\tilde n$: &         \num{1000} &  \num{10000} & \num{100000} & \num{1000000} & \num{10000000}\tabularnewline
			\midrule
			Sinkhorn's Algorithm~\ref{alg:Sinkhorn__sin_div}
						& 0.39\,s & \num{22.9}\,s & \multicolumn{3}{c}{\emph{out of memory or > 1\,hour}} \tabularnewline
			NFFT-accelerated Sinkhorn~\ref{alg:SinkhornFFT} 
						& 0.31\,s & \num{0.42}\,s & \num{2.0}\,s & \num{7.2}\,s & 59.4\,s \tabularnewline[2\doublerulesep]
			\bottomrule
		\end{tabular}
	\par\end{centering}

\smallskip

	\caption{\label{tab:Sinkhorn2}Two dimensions: comparison of computation times for $r=2$ and $\lambda=20$}
\end{table}

\subsection{Benchmark datasets}

\begin{figure}[htb!]
	\begin{center}
		
		{\includegraphics[width=0.6\textwidth, height=0.6\textwidth]{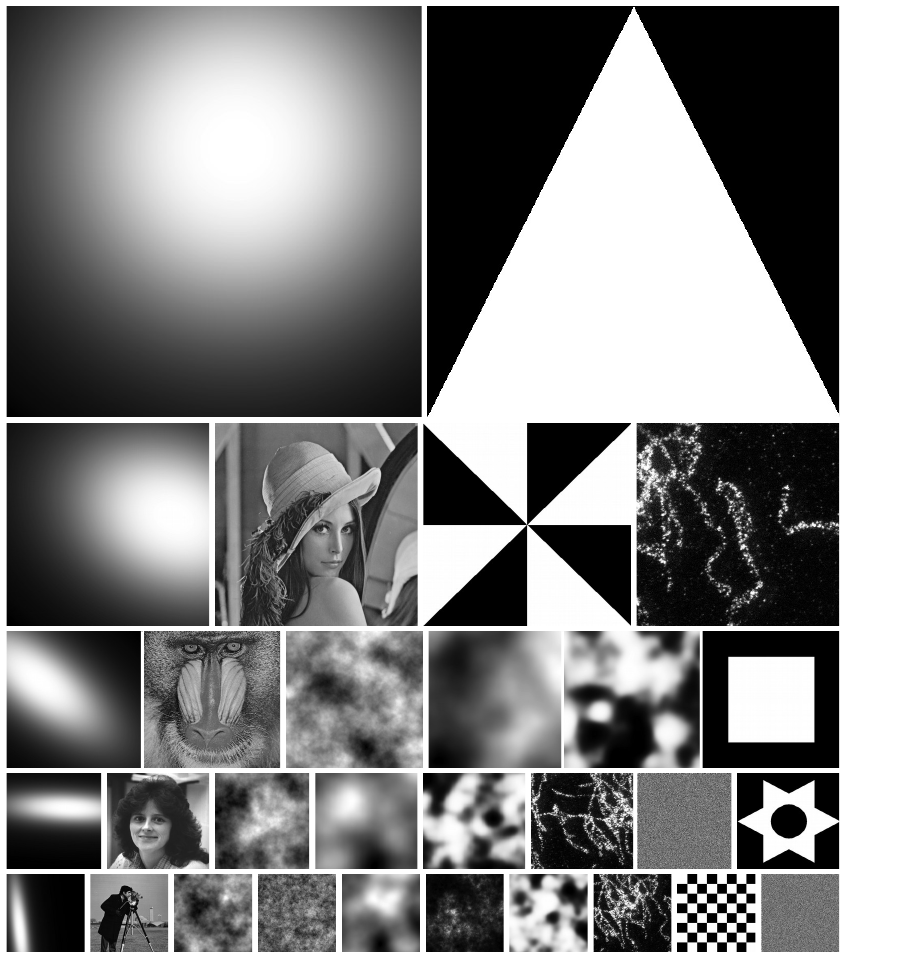}}
		
			\caption{Examples of the DOTmark database (from low to high resolution images) }\label{Fig:DOT Example}
	\end{center}
\end{figure}

This section validates the regularization parameter $\lambda$, and demonstrates the performance and the accuracy of NFFT-accelerated Sinkhorn's Algorithm~\ref{alg:SinkhornFFT} using real datasets.
We use a dataset called DOTmark (see Figure~\ref{Fig:DOT Example}); DOT stands for \emph{discrete optimal transport}.
This benchmark dataset is specially designed to effectively test and compare optimal transport methods~(cf.\ \citet{schrieber2016dotmark}).
It has gray level representation of the images in the resolution of $32\times 32$ to $512\times 512$, and it consists of 10 subsets of dataset, ranging from smooth to rough structure.

\paragraph{Transformation of images to probability vectors.}
\label{trans_img_prob}
\begin{figure}[htb!]
	\begin{center}
		\includegraphics[width=12cm ]{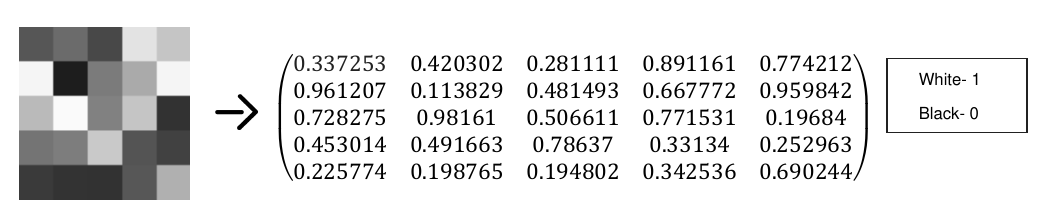}
		\caption{Grayscale image ($5\times 5$) is represented as a matrix}\label{Fig:Image to array}
	\end{center}
\end{figure}
A grayscale digital image can be represented as a matrix, where each entry represents a pixel in the image and the value of the pixel is the image's gray scale level in the range $[0,1]$ (see Figure~\ref{Fig:Image to array}).
In order to convert the grayscale image matrices into probability vectors, we vectorize and normalize the matrices. Furthermore, intensities of background pixels are the $\ell_1$ distance between pixels~$i$ and~$j$ of the respective grids ($32 \times 32,\dots, 512 \times 512$).

%
%
\newpage
\subsubsection{Validation of the regularization parameter \texorpdfstring{$\lambda$}{l}}
\label{Subsec:val_lamda}
\begin{figure}[htb!]
	\begin{center}
		
		{\includegraphics[width=1.\textwidth, height=0.7\textwidth]{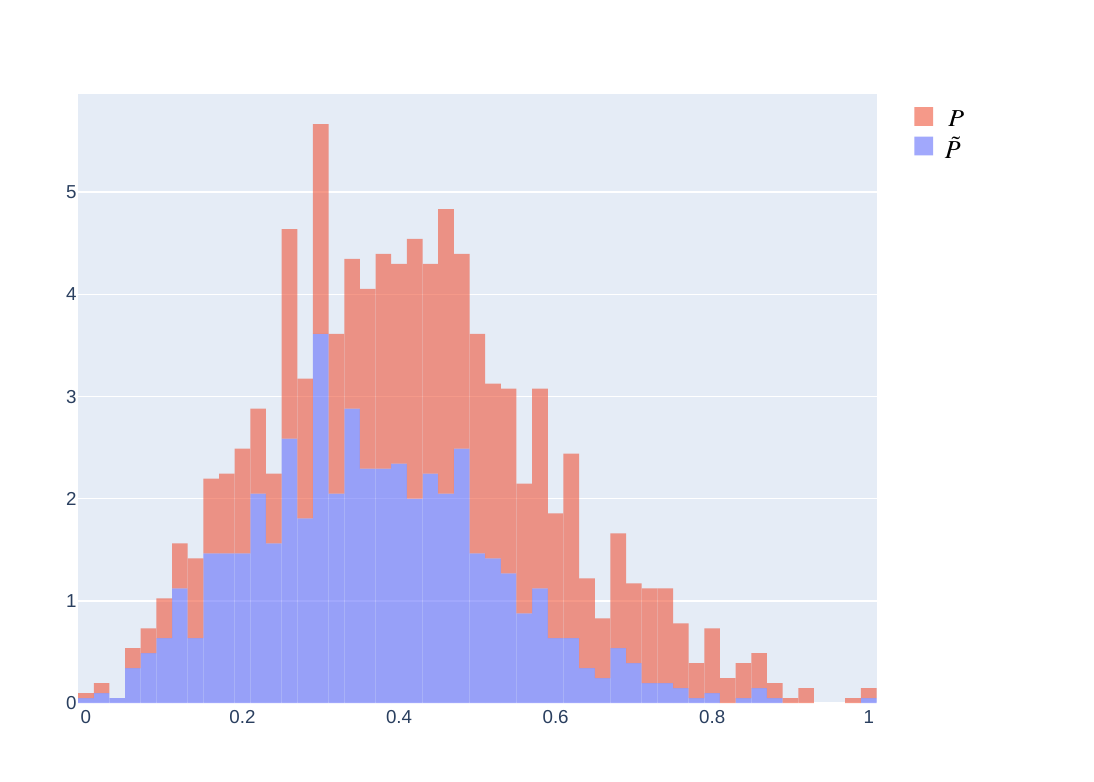}}
			\caption{Histogram of the measures~$P$ and~$\tilde P$ of the dataset GRFrough; $n$ ($\tilde n$, resp.) $= 1024$} \label{Fig:measures}
	\end{center}
\end{figure}
\begin{figure}[!htb]
	\centering
		\subfloat[Numerical results of Algorithm~\ref{alg:Sinkhorn__sin_div}]
		{\includegraphics[trim=0 0 0 10, clip, width=0.68\textwidth]{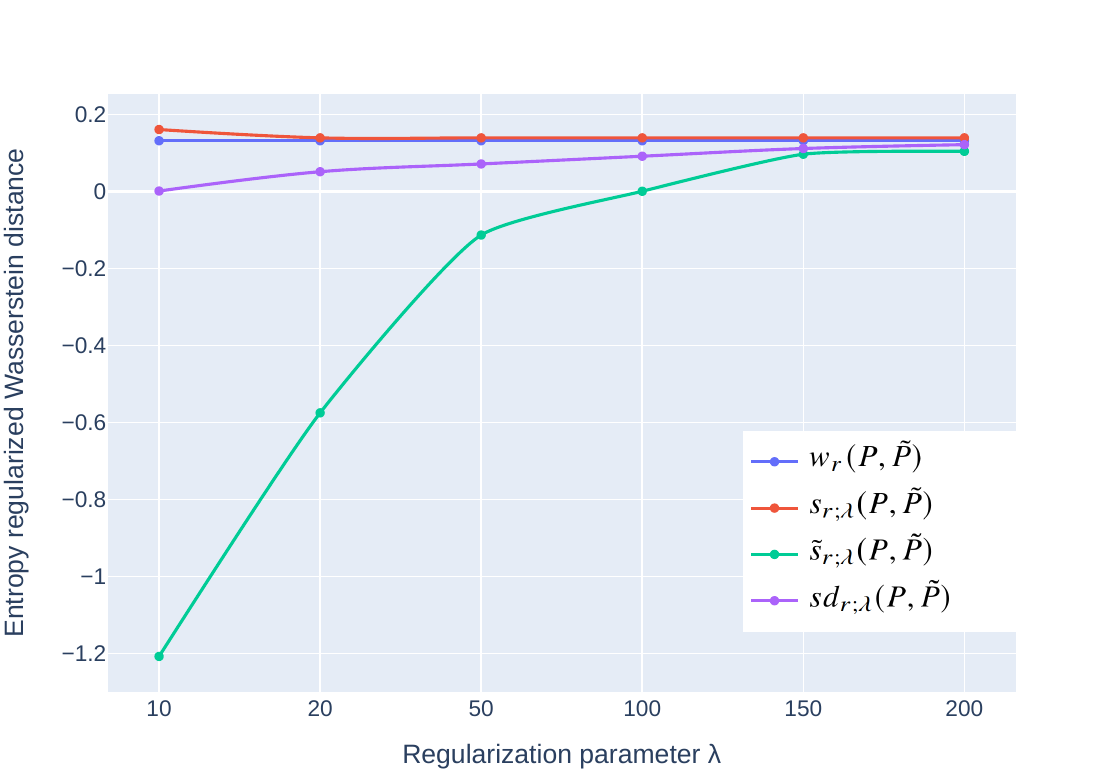} }
	\qquad
		\subfloat[Numerical results of Algorithm~\ref{alg:SinkhornFFT}]
		{\includegraphics[trim=0 0 0 10, clip, width=0.68\textwidth]{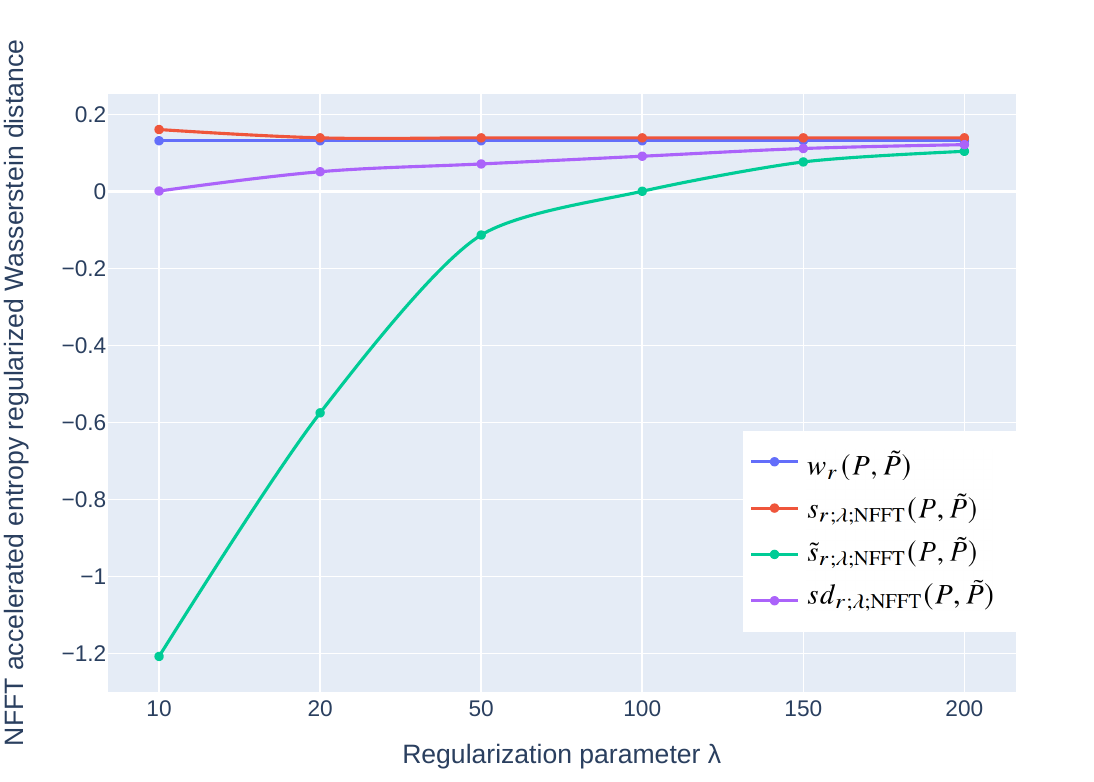}}
		\caption{Approximation of the Wasserstein distance using $s_{r;\lambda}(P,\tilde P)$, $\tilde{s}_{r;\lambda}(P,\tilde P)$, ${sd}_{r;\lambda}(P,\tilde P)$ and $s_{r;\lambda;\text{NFFT}}(P,\tilde P)$, $\tilde{s}_{r;\lambda;\text{NFFT}}(P,\tilde P)$, ${sd}_{r;\lambda;\text{NFFT}}(P,\tilde P)$. The parameters are $\lambda\in\{10,20,50,\dots,200\}$ and $r=2$.\label{Fig:SinkhornUB}}
\end{figure}
In this section we capture the behaviors of the lower and upper bounds (cf.\ Lemma~\ref{lem:Sinkhorn Approximation}) and Sinkhorn divergence, i.e., $s_{r;\lambda}(P,\tilde P)$, $\tilde s_{r;\lambda}(P,\tilde P)$, ${sd}_{r;\lambda}(P,\tilde P)$ and $s_{r;\lambda;\text{NFFT}}(P,\tilde P)$, $\tilde{s}_{r;\lambda;\text{NFFT}}(P,\tilde P)$, ${sd}_{r;\lambda;\text{NFFT}}(P,\tilde P)$ for different
values of the entropy regularization parameter $\lambda \in\{10,20,50,\dots,200\}$.
\newpage
We use the ‘GRFrough’ dataset, which is a subset of images from the DOTmark dataset. Notably, it has a rough structure, relative to the other subset of images (see Figure~\ref{Fig:measures}).
Figure~\ref{Fig:SinkhornUB} below investigate these quantities with respect to increasing~$\lambda$.
We infer that $s_{r;\lambda}(P,\tilde P)$, $s_{r;\lambda;\text{NFFT}}(P,\tilde P)$ converge slowly, for $\lambda$ increasing, to $w_r(P,\tilde P).$
However, $\tilde{s}_{r;\lambda}(P,\tilde P)$ and $\tilde{s}_{r;\lambda;\text{NFFT}}(P,\tilde P)$ converge quickly to $w_r(P,\tilde P),$ and Sinkhorn divergence $\big({sd}_{r;\lambda}(P,\tilde P)$, ${sd}_{r;\lambda;\text{NFFT}}(P,\tilde P)\big)$ also converge quicker, in comparison to $s_{r;\lambda}(P,\tilde P)$, $s_{r;\lambda;\text{NFFT}}(P,\tilde P).$ 
The argument behind these behaviors is, for larger $\lambda$, the weightage of the entropy in the objective function~\eqref{eq:Sinkhorn} decreases, and the matrices $\pi^s$ and $\pi^w$ coincide.
Furthermore, the NFFT approximation is stable for different values of regularization parameter~$\lambda$.
\bigskip


\subsubsection{Performance analysis}
\label{subsec:Performance analysis}
This section extensively substantiates the performance of NFFT-accelerated Sinkhorn's Algorithm~\ref{alg:SinkhornFFT} in terms of time and memory allocation.
For the experiments, we use the DOTmark dataset, ranging from 32 $\times$ 32 to 512 $\times$ 512 pixels in size, and  we consider transports between two different images of equal size.
\begin{table}[!htb]
	\footnotesize
	\begin{centering}
		\begin{tabular}{lcccccccccc}
			\toprule
			$n=\tilde n$ & \multicolumn{2}{c}{\num{1024}} & \multicolumn{2}{c}{\num{4096}} & \multicolumn{2}{c}{\num{16384}}& \multicolumn{2}{c}{\num{65536}} & \multicolumn{2}{c}{\num{262144}} \tabularnewline
			\midrule
			Dataset: \textbf{DOTmark} &  Alg.~\ref{alg:Sinkhorn__sin_div} & Alg.~\ref{alg:SinkhornFFT}&  Alg.~\ref{alg:Sinkhorn__sin_div} & Alg.~\ref{alg:SinkhornFFT}&  Alg.~\ref{alg:Sinkhorn__sin_div} & Alg.~\ref{alg:SinkhornFFT}&  Alg.~\ref{alg:Sinkhorn__sin_div} & Alg.~\ref{alg:SinkhornFFT}&  Alg.~\ref{alg:Sinkhorn__sin_div} & Alg.~\ref{alg:SinkhornFFT}  \tabularnewline
			\midrule
			CauchyDensity
			& \num{0.79}\,s & \textbf{\num{0.32}}\,s & \num{3.91}\,s & \textbf{\num{0.34}}\,s& \num{75.0}\,s & \textbf{\num{1.08}}\,s& - & \textbf{\num{1.11}}\,s& - & \textbf{\num{3.70}}\,s
			\tabularnewline
			ClassicImages
			& \num{0.75}\,s & \textbf{\num{0.29}}\,s & \num{3.53}\,s & \textbf{\num{0.36}}\,s& \num{69.0}\,s & \textbf{\num{1.11}}\,s& - & \textbf{\num{1.34}}\,s& - & \textbf{\num{3.69}}\,s
			\tabularnewline
			GRFmoderate
			& \num{1.51}\,s & \textbf{\num{0.41}}\,s & \num{4.15}\,s & \textbf{\num{0.53}}\,s& \num{93.16}\,s & \textbf{\num{1.26}}\,s& - & \textbf{\num{1.11}}\,s& - & \textbf{\num{3.36}}\,s
			\tabularnewline
			GRFrough
			& \num{0.89}\,s & \textbf{\num{0.39}}\,s & \num{3.72}\,s & \textbf{\num{0.46}}\,s& \num{79.0}\,s & \textbf{\num{1.31}}\,s& - & \textbf{\num{1.64}}\,s& - & \textbf{\num{3.81}}\,s
			\tabularnewline
			GRFsmooth
			& \num{1.92}\,s & \textbf{\num{0.43}}\,s & \num{5.15}\,s & \textbf{\num{0.61}}\,s& \num{102.13}\,s & \textbf{\num{1.46}}\,s& - & \textbf{\num{1.81}}\,s& - & \textbf{\num{3.52}}\,s
			\tabularnewline
			LogGRF
			& \num{2.19}\,s & \textbf{\num{0.52}}\,s & \num{5.31}\,s & \textbf{\num{0.69}}\,s& \num{105.10}\,s & \textbf{\num{1.78}}\,s& - & \textbf{\num{2.11}}\,s& - & \textbf{\num{4.71}}\,s
			\tabularnewline
			LogitGRF
			& \num{2.01}\,s & \textbf{\num{0.49}}\,s & \num{5.22}\,s & \textbf{\num{0.63}}\,s& \num{104.20}\,s & \textbf{\num{1.58}}\,s& - & \textbf{\num{2.01}}\,s& - & \textbf{\num{4.62}}\,s
			\tabularnewline
			MicroscopyImages
			& \num{0.53}\,s & \textbf{\num{0.21}}\,s & \num{3.34}\,s & \textbf{\num{0.23}}\,s& \num{67.9}\,s & \textbf{\num{0.72}}\,s& - & \textbf{\num{0.91}}\,s& - & \textbf{\num{1.86}}\,s
			\tabularnewline
			Shapes
			& \num{0.61}\,s & \textbf{\num{0.26}}\,s & \num{3.64}\,s & \textbf{\num{0.29}}\,s& \num{72.4}\,s & \textbf{\num{0.92}}\,s& - & \textbf{\num{1.01}}\,s& - & \textbf{\num{2.56}}\,s
			\tabularnewline
			WhiteNoise
			& \num{0.63}\,s & \textbf{\num{0.29}}\,s & \num{3.82}\,s & \textbf{\num{0.31}}\,s& \num{73.0}\,s & \textbf{\num{1.02}}\,s& - & \textbf{\num{1.07}}\,s& - & \textbf{\num{2.70}}\,s
			\tabularnewline
			\bottomrule
		\end{tabular}
	\par\end{centering}
\smallskip
	\caption{Comparison of computational time allocation of Sinkhorn's Algorithm~\ref{alg:Sinkhorn__sin_div} and NFFT-accelerated Sinkhorn's Algorithm~\ref{alg:SinkhornFFT}; $\lambda = 20$ and~$r=2$}  \label{tab:time_complexity}
\end{table}

\begin{table}[!htb]
	\footnotesize
	\begin{centering}
		\begin{tabular}{lcccccccccc}
			\toprule
			$n=\tilde n$ & \multicolumn{2}{c}{\num{1024}} & \multicolumn{2}{c}{\num{4096}} & \multicolumn{2}{c}{\num{16384}}& \multicolumn{2}{c}{\num{65536}} & \multicolumn{2}{c}{\num{262144}} \tabularnewline
			\midrule
			Dataset &  Alg.~\ref{alg:Sinkhorn__sin_div} & Alg.~\ref{alg:SinkhornFFT}&  Alg.~\ref{alg:Sinkhorn__sin_div} & Alg.~\ref{alg:SinkhornFFT}&  Alg.~\ref{alg:Sinkhorn__sin_div} & Alg.~\ref{alg:SinkhornFFT}&  Alg.~\ref{alg:Sinkhorn__sin_div} & Alg.~\ref{alg:SinkhornFFT}&  Alg.~\ref{alg:Sinkhorn__sin_div} & Alg.~\ref{alg:SinkhornFFT}  \tabularnewline
			\midrule
			 \textbf{DOTmark}& \multicolumn{2}{c}{(MB)} & \multicolumn{2}{c}{(MB)} & \multicolumn{2}{c}{(MB)}& \multicolumn{2}{c}{(MB)} & \multicolumn{2}{c}{(MB)} \tabularnewline
			\midrule
			CauchyDensity
			& \num{363.63} & \textbf{\num{37.11}} & \num{924.74} &\textbf{\num{51.94}} &\num{9032.14} &\textbf{\num{105.99}}& - & \textbf{\num{289.99 }}&- & \textbf{\num{321.23}} 
			\tabularnewline
			ClassicImages
			& \num{353.14} & \textbf{\num{35.31}} & \num{913.16} &\textbf{\num{49.14}} &\num{9000.11} &\textbf{\num{104.19}}& - & \textbf{\num{274.29 }}&- & \textbf{\num{311.29}} 
			\tabularnewline
			GRFmoderate
			& \num{412.12} & \textbf{\num{41.21}} & \num{941.23} &\textbf{\num{65.54}} &\num{9420.35} &\textbf{\num{135.17}}& - & \textbf{\num{301.15}}&- & \textbf{\num{333.52}} 
			\tabularnewline
			GRFrough
			& \num{341.34} & \textbf{\num{33.11}} & \num{912.13} &\textbf{\num{43.12}} &\num{8909.13} &\textbf{\num{99.12}}& - & \textbf{\num{263.21 }}&- & \textbf{\num{301.21}} 
			\tabularnewline
			GRFsmooth
			& \num{442.13} & \textbf{\num{51.21}} & \num{981.41} &\textbf{\num{85.34}} &\num{9740.31} &\textbf{\num{149.19}}& - & \textbf{\num{331.19}}&- & \textbf{\num{373.12}} 
			\tabularnewline
			LogGRF
			& \num{463.63} & \textbf{\num{57.11}} & \num{1124.14} &\textbf{\num{91.94}} &\num{9831.34} &\textbf{\num{155.99}}& - & \textbf{\num{359.19}}&- & \textbf{\num{391.95}} 
			\tabularnewline
			LogitGRF
			& \num{451.23} & \textbf{\num{52.41}} & \num{1023.11} &\textbf{\num{89.14}} &\num{9800.54} &\textbf{\num{151.49}}& - & \textbf{\num{349.29}}&- & \textbf{\num{384.13}} 
			\tabularnewline
			MicroscopyImages
			& \num{250.11} & \textbf{\num{19.41}} & \num{912.17} &\textbf{\num{35.17}} &\num{8611.39} &\textbf{\num{55.12}}& - & \textbf{\num{221.74 }}&- & \textbf{\num{243.21}} 
			\tabularnewline
			Shapes
			& \num{311.17} & \textbf{\num{29.41}} & \num{922.12} &\textbf{\num{41.13}} &\num{8751.27} &\textbf{\num{95.13}}& - & \textbf{\num{241.21 }}&- & \textbf{\num{283.61}} 
			\tabularnewline
			WhiteNoise
			& \num{321.12} & \textbf{\num{31.41}} & \num{932.13} &\textbf{\num{42.11}} &\num{8800.12} &\textbf{\num{97.11}}& - & \textbf{\num{253.22 }}&- & \textbf{\num{299.31}} 
			\tabularnewline
			\bottomrule
		\end{tabular}
	\par\end{centering}
\smallskip
	\caption{Comparison of computational memory allocation of Sinkhorn's Algorithm~\ref{alg:Sinkhorn__sin_div} and NFFT-accelerated Sinkhorn's Algorithm~\ref{alg:SinkhornFFT}; $\lambda = 20$ and~$r=2$  }\label{tab:memory_allocations}
\end{table}

Table~\ref{tab:time_complexity} and~\ref{tab:memory_allocations} compare the time and memory allocation of Sinkhorn's Algorithm~\ref{alg:Sinkhorn__sin_div} and NFFT accelarated Sinkhorn's Algorithm~\ref{alg:SinkhornFFT}. Among all the tests, the computational time and memory allocation of our proposed algorithm is significantly better. 

From the performance analysis, we infer that NFFT-accelerated Sinkhorn's Algorithm~\ref{alg:SinkhornFFT} utilizes lesser time and memory, in spite of high number of points $n$ ($\tilde{n},$ resp.)\ for computations, and our device runs out of memory for Sinkhorn's Algorithm~\ref{alg:Sinkhorn__sin_div}, when the problems are sized larger than 16384 $n$ ($\tilde{n},$ resp.).


\subsubsection{Accuracy analysis}\label{sec:Accuracy}
We validate the computational accuracy of NFFT-accelerated Sinkhorn's Algorithm~\ref{alg:SinkhornFFT}.
Throughout the accuracy analysis, we use $\tilde{s}_{r;\lambda}(P,\tilde P)$ and $\tilde{s}_{r;\lambda;\text{NFFT}}(P,\tilde P) $, since it is a better approximation of Wasserstein distance (see Section~\ref{Subsec:val_lamda}).
Initially, we perform the accuracy analysis using the low resolution images from the DOTmark dataset. 
From Table~\ref{Tab:Accuracy_analysis2}, we notice that NFFT-accelerated Sinkhorn's Algorithm~\ref{alg:SinkhornFFT} has achieved stable approximation without compromise in accuracy.
\begin{table}[!htb]
		\begin{centering}
		\begin{tabular}{lccc}
			\toprule
			$n\times\tilde n=1024\times 1024$ &  Wasserstein & Sinkhorn~\ref{alg:Sinkhorn__sin_div} & NFFT-accelerated Sinkhorn~\ref{alg:SinkhornFFT} \tabularnewline
			\midrule
			Dataset: \textbf{DOTmark}&  $w_r(P,\tilde P)$ &  $\tilde{s}_{r;\lambda}(P,\tilde P)$ & $\tilde{s}_{r;\lambda;\text{NFFT}}(P,\tilde P)$ \tabularnewline
			\midrule
			CauchyDensity
			& \num{0.1204981}& \num{0.1204992} & \num{0.1204992}\tabularnewline
			ClassicImages 
			& \num{0.0629847}& \num{0.0629848} & \num{0.0629848}\tabularnewline
			GRFmoderate
			& \num{0.0600861}& \num{0.0609624}&\num{0.0609624}\tabularnewline
			GRFrough
			& \num{0.0324286}& \num{0.0327149} &\num{0.0327149}\tabularnewline
			GRFsmooth
			& \num{0.1422346}& \num{0.1456118} & \num{0.1456118}
			\tabularnewline
			LogGRF	
			& \num{0.1260294} & \num{0.1267356} &\num{0.1267356} 
			\tabularnewline
			LogitGRF	
			& \num{0.1075698} & \num{0.1077784} &\num{0.1077784} 
			\tabularnewline
			MicroscopyImages	
			& \num{0.0920382} & \num{0.0924179} & \num{0.0924179}
			\tabularnewline
			Shapes	
			& \num{0.1435847} & \num{0.1447279} & \num{0.1447279}
			\tabularnewline
			WhiteNoise	
			& \num{0.0206093} & \num{0.0206900} & \num{0.0206900}
			\tabularnewline
			\bottomrule
		\end{tabular}
	\par\end{centering}

\smallskip

	\caption{Comparison of accuracy to compute or approximate the Wasserstein distance in low-resolution images; $\lambda = 20$ and~$r=2$ \label{Tab:Accuracy_analysis2}}
\end{table}
\begin{table}[!htb]
	\footnotesize
	\begin{centering}
		\begin{tabular}{lcccccc}
			\toprule
			Dataset: GRFrough & \multicolumn{2}{c}{Wasserstein} & \multicolumn{2}{c}{Sinkhorn~\ref{alg:Sinkhorn__sin_div}} & \multicolumn{2}{c}{NFFT-accelerated Sinkhorn~\ref{alg:SinkhornFFT}} \tabularnewline
			\midrule
			$n\times\tilde n$ &  $w_r(P,\tilde P)$ & time & $\tilde{s}_{r;\lambda}(P,\tilde P)$ & time  & $\tilde{s}_{r;\lambda;\text{NFFT}}(P,\tilde P)$ &time  \tabularnewline
			\midrule
			$1024\times1024$
			& \num{0.0324286}& 72.34\,s &\num{0.0327149} & 0.89\,s& \num{0.0327149}& 0.39\,s \tabularnewline
			$4096\times4096$ 
			& \multicolumn{2}{c}{\emph{out of memory}}&  \num{0.0872946} &3.72\,s &\num{0.0872946}&0.46\,s \tabularnewline
			$16384\times16384$
			& \multicolumn{2}{c}{\emph{out of memory}} &\num{0.1249751}& 79.0\,s&\num{0.1249751}&1.31\,s \tabularnewline
			$65536\times65536$
			& \multicolumn{2}{c}{\emph{out of memory}} &\multicolumn{2}{c}{\emph{out of memory}}&\num{0.7432192}&1.64\,s \tabularnewline
			$262144\times262144$
			& \multicolumn{2}{c}{\emph{out of memory}} &\multicolumn{2}{c}{\emph{out of memory}}& \num{0.9117436} & 3.81\,s  
			\tabularnewline
			\bottomrule
		\end{tabular}
	\par\end{centering}
\smallskip
	\caption{Comparison of computation times and accuracy to compute or approximate the Wasserstein distance, from low to high resolution images; $\lambda = 20$ and~$r=2$}  \label{Tab:Accuracy_analysis1}
\end{table}

Now, we move on to high resolution images of the ‘GRFrough’ dataset.
Table~\ref{Tab:Accuracy_analysis1} comprises the list of values that enables us to understand the approximation accuracy, as we move from low to high resolution images.
We compare the results of Sinkhorn's Algorithm~\ref{alg:Sinkhorn__sin_div} with NFFT-accelerated Sinkhorn's Algorithm~\ref{alg:SinkhornFFT} for the problems sized up to \num{16384} $n$ ($\tilde{n},$ resp.), and we infer that there is no compromise in accuracy.
We recognize that advancing to high resolution images does not affect the stability of approximation.
Due to the break of Sinkhorn's Algorithm~\ref{alg:Sinkhorn__sin_div}, it cannot be used as comparison factor when the size of the problem is beyond \num{16384} $n$ ($\tilde{n},$ resp.). However, our proposed algorithm computes the largest problem available in the DOTmark dataset, which is of size \num{262144}
 $n$ ($\tilde{n},$ resp.).



\subsection{Comparisons and further steps} \label{sec:hist_rem}

In this section, we substantiate the historical evolution of the prominent algorithms, which approximate the Wasserstein distance.
Furthermore, we discuss the supremacy and the direction of further development of our proposed algorithms. 
\paragraph{Historical remarks.} The approach of entropy regularization of the Wasserstein distance by \citet{cuturi2013sinkhorn} is a well-known path breaking approach to approximate the Wasserstein distance, which is effectively computed by Sinkhorn's algorithm.
Later on, many constructive approaches and/\,or analyses were contributed to improve and/\,or support the entropy regularization approach (cf.\ \citet{AltschulerSinkhorn}, \citet{dvurechensky2018computational}). 
In 2019, the approach of log-domain stabilization and truncated kernel of the Sinkhorn's algorithm was proposed by \citet{schmitzer2019stabilized}. The log-domain stabilization method satisfies the demand for larger regularization parameters~$\lambda$, and the truncated kernel reduces memory demand and also accelerates the iterations.
In the same article, a multi-scale scheme was also proposed, which enables more efficient computations of the kernel truncated approach.
As discussed in Remark~\ref{Sec:sinkhorn_div}, these prominent approaches still suffer by the entropy bias. In order to remove/\,reduce the bias, Sinkhorn divergence was proposed by \citet{ramdas2017wasserstein}.

These prominent approaches affirm the progressive improvement of the algorithm, which approximates the Wasserstein distance.
However, notably, these approaches still deteriorate by the matrix-vector operations, which is the bottleneck of the algorithms. 
\begin{table}[!htb]
	\begin{centering}
		\begin{tabular}{lllll}
			\toprule 
			Name of algorithm / method &         Algorithm &  denotement  \tabularnewline
			\midrule
			Standard Sinkhorn~(\citet{cuturi2013sinkhorn})
			& Algorithm~\ref{alg:Sinkhorn__sin_div} & Std. Sinkhorn \tabularnewline
			Stabilized log-domain Sinkhorn~(\citet{schmitzer2019stabilized}) 
			&Algorithm~\ref{alg:Sinkhorn_Log_sin} & Stb. log Sinkhorn \tabularnewline
			Sinkhorn divergence~(\citet{ramdas2017wasserstein})
			& Algorithm~\ref{alg:Sinkhorn__sin_div} & ${sd}_{r;\lambda}(P,\tilde P)$
			\tabularnewline
			Multi-scale Sinkhorn~(\citet{schmitzer2019stabilized})
			& Algorithm  \footnote[9]{ Cf.\ \href{https://github.com/OTGroupGoe/MultiScaleOT.jl}{https://github.com/OTGroupGoe/MultiScaleOT.jl}}  & Multi Sinkhorn   
			\tabularnewline
			\midrule
			NFFT-accelerated Sinkhorn
			& Algorithm~\ref{alg:SinkhornFFT} & NFFT Sinkhorn   
			\tabularnewline
			NFFT-accelerated log-domain Sinkhorn
			& Algorithm~\ref{alg:logSinkhornFFT} & NFFT log Sinkhorn   
			\tabularnewline
			NFFT-accelerated Sinkhorn divergence
			& Algorithm~\ref{alg:SinkhornFFT} & ${sd}_{r;\lambda;\text{NFFT}}(P,\tilde P)$  
			\tabularnewline
			\bottomrule
		\end{tabular}
	\end{centering}

\smallskip
	\caption{List of algorithms\label{Tab:ListofAlg}  }
\end{table}

\begin{figure}[!htb]
	\begin{center}
		\includegraphics[width=15cm]{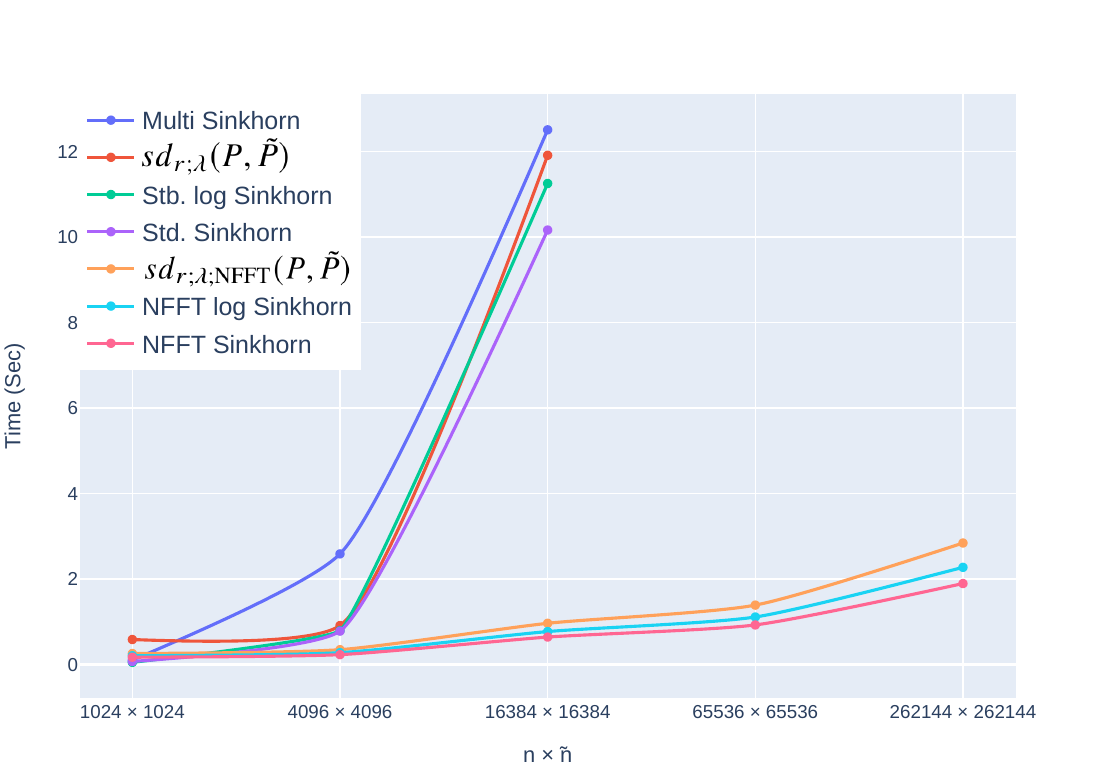}
		\caption{Comparison of computational time allocations of algorithms in Table~\ref{Tab:ListofAlg}; $\lambda = 20$ and~$r=2$}\label{Fig:Time_analysis_diff_alg}
		\end{center}
\end{figure}

\begin{table}[!htb]
	\begin{centering}
		\begin{tabular}{lc c c c c}
			\toprule 
			$n=\tilde n$ &\num{1024}&\num{4096}&\num{16384}&\num{65536}&\num{262144}  \tabularnewline
			\midrule
			Algorithm & (MB)& (MB)& (MB)& (MB)& (MB)\tabularnewline
			\midrule
			Std. Sinkhorn
			& \num{27.44} & \num{381.10} & \num{6644.10}& \multicolumn{2}{c}{\emph{out of memory}}\tabularnewline
			Stb. log Sinkhorn
			& \num{103.10} & \num{487.48} & \num{6943.12} & \multicolumn{2}{c}{\emph{out of memory}} \tabularnewline
			${sd}_{r;\lambda}(P,\tilde P)$
			& \num{51.34} & \num{391.49}& \num{7139.12} & \multicolumn{2}{c}{\emph{out of memory}}  \tabularnewline
			Multi Sinkhorn 
			&\num{1.9} &\num{19.49} &\num{32.13} & \multicolumn{2}{c}{\emph{out of memory}}   
			\tabularnewline
			\midrule
			NFFT Sinkhorn
			& \textbf{\num{1.7}} &\textbf{\num{2.06}}& \textbf{\num{8.98}} &\textbf{\num{37.5}} & \textbf{\num{132.0}}   \tabularnewline 
			NFFT log Sinkhorn
			& \textbf{\num{1.79}} &\textbf{\num{2.10}}& \textbf{\num{9.31}} &\textbf{\num{43.1}} & \textbf{\num{142.7}}      \tabularnewline 
			${sd}_{r;\lambda;\text{NFFT}}(P,\tilde P)$
			& \textbf{\num{2.3}} &\textbf{\num{3.19}}& \textbf{\num{12.78}} &\textbf{\num{45.29}} & \textbf{\num{162.0}}   
		\tabularnewline 
			\bottomrule
		\end{tabular}
	\par\end{centering}

\smallskip

	\caption{Comparison of computational memory allocations of algorithms in Table~\ref{Tab:ListofAlg}; $\lambda = 20$ and~$r=2$  \label{Tab:memory_analysis_diff_alg}}
\end{table}
Now, we compare our proposed algorithms with prominent algorithms, which are discussed so far (see Table~\ref{Tab:ListofAlg}). We would like to emphasize that our proposed algorithms are compatible even with low-threshold applications, this does
not require expensive hardware or having access to supercomputers.
All the algorithms involved in the comparison including our proposed algorithms follow Central Processing Unit~(CPU) implementation paradigms. 
We follow the same experimentally setup utilized in preceding Section~\ref{sec:Accuracy}, and we use the 'GRFrough' dataset.
From Figure~\ref{Fig:Time_analysis_diff_alg} and Table~\ref{Tab:memory_analysis_diff_alg}, it is evident that our proposed algorithms perform significantly better in terms of time and memory allocations. 
Our device runs out of memory for all algorithms/\,methods, except our proposed algorithms, when the problems are sized larger than \num{16384} $n$ ($\tilde{n},$ resp.).
In terms of memory allocations, the Multi Sinkhorn algorithm shows significant performance, and the results are closer to our proposed algorithms.
However, it requires more computational time, and it breaks due to the kernel matrix formation, when the problems are sized larger than \num{16384} $n$ ($\tilde{n},$ resp.).

\paragraph{Faster computation.} In general, for faster computations, Graphics Processing Unit~(GPU) implementations are used.
The 'GeomLoss' package is a clever GPU implementation to approximate the  Wasserstein distance or to solve the OT problem.
We refer to \citet{feydy2020geometric} and the corresponding GitHub repository for further information on the implementation.
As mentioned in Section~\ref{sec:intro}, for a fast computation, low rank approximation techniques are also considered. 
However, the 'GeomLoss' is the prominent contribution in terms of fast computation.
\paragraph{Further steps.} Our proposed algorithms surpass the burden of time and memory allocations, and it is also flexible to adapt to the log-domain implementation.
However, further research directions will be focused on applying our algorithms to suitable applications, and incorporation of possible extensions.
In order to reach wider audiences, the GPU implementation of our proposed algorithms can be considered as one of the possible extensions as well. We would like to mention that GPU implementation of NFFT algorithm is readily available in  corresponding GitHub repository.\footnote{Cf.\ \href{https://github.com/sukunis/CUNFFT/tree/master/src}{https://github.com/sukunis/CUNFFT/tree/master/src}}

\section{Summary} \label{Sec:summary}
The nonequispaced fast Fourier transform, as presented in this article, allows computing a proper approximation of the Wasserstein distance in $\mathcal O(n \log n )$ arithmetic operations.
NFFT-accelerated Sinkhorn's Algorithm~\ref{alg:SinkhornFFT} performs significantly better than standard Sinkhorn's Algorithm~\ref{alg:Sinkhorn__sin_div}, in terms of computational time and memory allocations. Our numerical results demonstrate the effectiveness of the new method as well as the tightness of our theoretical bounds. We believe that our algorithms can be widely used in data sciene applications for handling large-scale dataset.





\bibliographystyle{abbrvnat}
\bibliography{LiteraturAlois,LiteraturDaniel,LiteraturRaj}

\begin{thebibliography}{40}
\providecommand{\natexlab}[1]{#1}
\providecommand{\url}[1]{\texttt{#1}}
\expandafter\ifx\csname urlstyle\endcsname\relax
  \providecommand{\doi}[1]{doi: #1}\else
  \providecommand{\doi}{doi: \begingroup \urlstyle{rm}\Url}\fi

\bibitem[Altschuler et~al.(2017)Altschuler, Weed, and
  Rigollet]{AltschulerSinkhorn}
J.~Altschuler, J.~Weed, and P.~Rigollet.
\newblock Near-linear time approximation algorithms for optimal transport via
  {S}inkhorn iteration.
\newblock In \emph{Proceedings of the 31st International Conference on Neural
  Information Processing Systems}, pages 1961--1971. Curran Associates Inc.,
  2017.
\newblock URL \url{https://arxiv.org/abs/1705.09634}.

\bibitem[Altschuler et~al.(2019)Altschuler, Bach, Rudi, and
  Niles-Weed]{NEURIPS2019_f55cadb9}
J.~Altschuler, F.~Bach, A.~Rudi, and J.~Niles-Weed.
\newblock Massively scalable sinkhorn distances via the nystr\"{o}m method.
\newblock In H.~Wallach, H.~Larochelle, A.~Beygelzimer, F.~d\textquotesingle
  Alch\'{e}-Buc, E.~Fox, and R.~Garnett, editors, \emph{Advances in Neural
  Information Processing Systems}, volume~32. Curran Associates, Inc., 2019.
\newblock URL
  \url{https://proceedings.neurips.cc/paper/2019/file/f55cadb97eaff2ba1980e001b0bd9842-Paper.pdf}.

\bibitem[Altschuler and Boix-Adsera(2020)]{altschuler2020polynomial}
J.~M. Altschuler and E.~Boix-Adsera.
\newblock Polynomial-time algorithms for multimarginal optimal transport
  problems with structure.
\newblock \emph{arXiv preprint arXiv:2008.03006}, 2020.
\newblock URL \url{https://arxiv.org/abs/2008.03006}.

\bibitem[Ba and Quellmalz(2022)]{QuellmalzFFT}
F.~A. Ba and M.~Quellmalz.
\newblock Accelerating the {S}inkhorn algorithm for sparse multi-marginal
  optimal transport by fast {F}ourier transforms, 2022.
\newblock URL \url{https://arXiv.org/abs/2208.03120}.

\bibitem[Bilik et~al.(2019)Bilik, Longman, Villeval, and
  Tabrikian]{bilik2019rise}
I.~Bilik, O.~Longman, S.~Villeval, and J.~Tabrikian.
\newblock The rise of radar for autonomous vehicles: Signal processing
  solutions and future research directions.
\newblock \emph{IEEE signal processing Magazine}, 36\penalty0 (5):\penalty0
  20--31, 2019.
\newblock \doi{10.1109/MSP.2019.2926573}.

\bibitem[Bolley(2008)]{Bolley2008}
F.~Bolley.
\newblock Separability and completeness for the {W}asserstein distance.
\newblock In C.~Donati-Martin, M.~\'Emery, A.~Rouault, and C.~Stricker,
  editors, \emph{S\'eminaire de Probabilit\'es XLI}, volume 1934 of
  \emph{Lecture Notes in Mathematics}, pages 371--377. Springer, Berlin,
  Heidelberg, 2008.
\newblock \doi{10.1007/978-3-540-77913-1}.

\bibitem[Chakraborty et~al.(2020)Chakraborty, Paul, and
  Das]{chakraborty2020hierarchical}
S.~Chakraborty, D.~Paul, and S.~Das.
\newblock Hierarchical clustering with optimal transport.
\newblock \emph{Statistics \& Probability Letters}, 163:\penalty0 108781, 2020.
\newblock \doi{10.1016/j.spl.2020.108781}.

\bibitem[Courty et~al.(2014)Courty, Flamary, and Tuia]{courty2014domain}
N.~Courty, R.~Flamary, and D.~Tuia.
\newblock Domain adaptation with regularized optimal transport.
\newblock In \emph{Joint European Conference on Machine Learning and Knowledge
  Discovery in Databases}, pages 274--289. Springer, 2014.
\newblock \doi{10.1007/978-3-662-44848-9}.

\bibitem[Cuturi(2013)]{cuturi2013sinkhorn}
M.~Cuturi.
\newblock Sinkhorn distances: Lightspeed computation of optimal transport.
\newblock \emph{Advances in neural information processing systems}, 26, 2013.
\newblock URL \url{https://proceedings.mlr.press/v89/feydy19a.html}.

\bibitem[Dvurechensky et~al.(2018)Dvurechensky, Gasnikov, and
  Kroshnin]{dvurechensky2018computational}
P.~Dvurechensky, A.~Gasnikov, and A.~Kroshnin.
\newblock Computational optimal transport: Complexity by accelerated gradient
  descent is better than by {S}inkhorn’s algorithm.
\newblock In \emph{International conference on machine learning}, pages
  1367--1376. PMLR, 2018.
\newblock URL \url{http://proceedings.mlr.press/v80/dvurechensky18a.html}.

\bibitem[Feydy(2020)]{feydy2020geometric}
J.~Feydy.
\newblock \emph{Geometric data analysis, beyond convolutions}.
\newblock PhD thesis, Universit{\'e} Paris-Saclay Gif-sur-Yvette, France, 2020.
\newblock URL
  \url{https://www.math.ens.psl.eu/~feydy/geometric_data_analysis_draft.pdf}.

\bibitem[Gasnikov et~al.(2016)Gasnikov, Gasnikova, Nesterov, and
  Chernov]{gasnikov2016efficient}
A.~V. Gasnikov, E.~Gasnikova, Y.~E. Nesterov, and A.~Chernov.
\newblock Efficient numerical methods for entropy-linear programming problems.
\newblock \emph{Computational Mathematics and Mathematical Physics},
  56\penalty0 (4):\penalty0 514--524, 2016.
\newblock URL
  \url{https://link.springer.com/content/pdf/10.1134/S0965542516040084.pdf}.

\bibitem[Genevay(2019)]{genevay2019entropy}
A.~Genevay.
\newblock \emph{Entropy-regularized optimal transport for machine learning}.
\newblock PhD thesis, Paris Sciences et Lettres (ComUE), 2019.
\newblock URL \url{https://www.theses.fr/2019PSLED002}.

\bibitem[Graf and Luschgy(2000)]{GrafLuschgy}
S.~Graf and H.~Luschgy.
\newblock \emph{Foundations of Quantization for Probability Distributions},
  volume 1730 of \emph{Lecture Notes in Mathematics}.
\newblock Springer-Verlag, Berlin, 2000.
\newblock \doi{10.1007/BFb0103945}.

\bibitem[Hao et~al.(2013)Hao, Kilmer, Braman, and Hoover]{hao2013facial}
N.~Hao, M.~E. Kilmer, K.~Braman, and R.~C. Hoover.
\newblock Facial recognition using tensor-tensor decompositions.
\newblock \emph{SIAM Journal on Imaging Sciences}, 6\penalty0 (1):\penalty0
  437--463, 2013.
\newblock \doi{10.1137/110842570}.
\newblock URL \url{https://doi.org/10.1137/110842570}.

\bibitem[Kalantari et~al.(2008)Kalantari, Lari, Ricca, and
  Simeone]{kalantari2008complexity}
B.~Kalantari, I.~Lari, F.~Ricca, and B.~Simeone.
\newblock On the complexity of general matrix scaling and entropy minimization
  via the ras algorithm.
\newblock \emph{Mathematical Programming}, 112\penalty0 (2):\penalty0 371--401,
  2008.
\newblock \doi{10.1007/s10107-006-0021-4}.

\bibitem[Keiner et~al.()Keiner, Kunis, and Potts]{nfft3}
J.~Keiner, S.~Kunis, and D.~Potts.
\newblock {NFFT 3.5, C subroutine library}.
\newblock \url{http://www.tu-chemnitz.de/~potts/nfft}.
\newblock Contributors: F.~Bartel, M.~Fenn, T.~G\"orner, M.~Kircheis, T.~Knopp,
  M.~Quellmalz, M.~Schmischke, T.~Volkmer, A.~Vollrath.

\bibitem[Khalil~Abid and Gower(2018)]{khalil2018greedy}
B.~Khalil~Abid and R.~M. Gower.
\newblock Greedy stochastic algorithms for entropy-regularized optimal
  transport problems.
\newblock \emph{arXiv e-prints}, pages arXiv--1803, 2018.
\newblock URL \url{http://proceedings.mlr.press/v84/abid18a/abid18a.pdf}.

\bibitem[Kusner et~al.(2015)Kusner, Sun, Kolkin, and
  Weinberger]{kusner2015word}
M.~Kusner, Y.~Sun, N.~Kolkin, and K.~Weinberger.
\newblock From word embeddings to document distances.
\newblock In \emph{International conference on machine learning}, pages
  957--966. PMLR, 2015.
\newblock URL \url{https://proceedings.mlr.press/v37/kusnerb15.html}.

\bibitem[Lin et~al.(2019)Lin, Ho, and Jordan]{lin2019efficiency}
T.~Lin, N.~Ho, and M.~I. Jordan.
\newblock On the efficiency of {S}inkhorn and {G}reenkhorn and their
  acceleration for optimal transport.
\newblock \emph{arXiv preprint arXiv:1906.01437}, 2019.
\newblock URL \url{https://proceedings.mlr.press/v97/lin19a.html}.

\bibitem[Luise et~al.(2018)Luise, Rudi, Pontil, and
  Ciliberto]{NEURIPS2018_3fc2c60b}
G.~Luise, A.~Rudi, M.~Pontil, and C.~Ciliberto.
\newblock Differential properties of sinkhorn approximation for learning with
  wasserstein distance.
\newblock In S.~Bengio, H.~Wallach, H.~Larochelle, K.~Grauman, N.~Cesa-Bianchi,
  and R.~Garnett, editors, \emph{Advances in Neural Information Processing
  Systems}, volume~31. Curran Associates, Inc., 2018.
\newblock URL
  \url{https://proceedings.neurips.cc/paper/2018/file/3fc2c60b5782f641f76bcefc39fb2392-Paper.pdf}.

\bibitem[Mensch and Peyr{\'e}(2020)]{mensch2020online}
A.~Mensch and G.~Peyr{\'e}.
\newblock Online sinkhorn: Optimal transport distances from sample streams.
\newblock \emph{Advances in Neural Information Processing Systems},
  33:\penalty0 1657--1667, 2020.
\newblock URL
  \url{https://proceedings.neurips.cc/paper/2020/hash/123650dd0560587918b3d771cf0c0171-Abstract.html}.

\bibitem[Nestler(2018)]{nestlerdiss}
F.~Nestler.
\newblock \emph{{Efficient Computation of Electrostatic Interactions in
  Particle Systems Based on Nonequispaced Fast Fourier Transforms}}.
\newblock Dissertation. Universit\"atsverlag Chemnitz, 2018.
\newblock ISBN 978-3-96100-054-8.
\newblock URL \url{http://nbn-resolving.de/urn:nbn:de:bsz:ch1-qucosa2-233760}.

\bibitem[Neumayer and Steidl(2021)]{neumayer2021optimal}
S.~Neumayer and G.~Steidl.
\newblock From optimal transport to discrepancy.
\newblock \emph{Handbook of Mathematical Models and Algorithms in Computer
  Vision and Imaging: Mathematical Imaging and Vision}, pages 1--36, 2021.
\newblock URL
  \url{https://link.springer.com/content/pdf/10.1007/978-3-030-03009-4_95-1.pdf}.

\bibitem[Papadakis et~al.(2014)Papadakis, Peyr\'{e}, and
  Oudet]{papadakis2014optimal}
N.~Papadakis, G.~Peyr\'{e}, and E.~Oudet.
\newblock Optimal transport with proximal splitting.
\newblock \emph{SIAM Journal on Imaging Sciences}, 7\penalty0 (1):\penalty0
  212--238, 2014.
\newblock \doi{10.1137/130920058}.

\bibitem[Peyr{\'e} and Cuturi(2019)]{peyre2019computational}
G.~Peyr{\'e} and M.~Cuturi.
\newblock Computational optimal transport: With applications to data science.
\newblock \emph{Foundations and Trends{\textregistered} in Machine Learning},
  11\penalty0 (5--6):\penalty0 355--607, 2019.
\newblock URL \url{https://ieeexplore.ieee.org/document/8641476}.

\bibitem[Platte et~al.(2011)Platte, Trefethen, and
  Kuijlaars]{platte2011impossibility}
R.~B. Platte, L.~N. Trefethen, and A.~B. Kuijlaars.
\newblock Impossibility of fast stable approximation of analytic functions from
  equispaced samples.
\newblock \emph{SIAM review}, 53\penalty0 (2):\penalty0 308--318, 2011.
\newblock \doi{10.1137/090774707}.

\bibitem[Plonka et~al.(2018)Plonka, Potts, Steidl, and Tasche]{PlPoStTa18}
G.~Plonka, D.~Potts, G.~Steidl, and M.~Tasche.
\newblock \emph{Numerical Fourier Analysis}.
\newblock Applied and Numerical Harmonic Analysis. Birkh\"auser, 2018.
\newblock ISBN 978-3-030-04305-6.
\newblock \doi{10.1007/978-3-030-04306-3}.

\bibitem[Ramdas et~al.(2017)Ramdas, Garc{\'\i}a~Trillos, and
  Cuturi]{ramdas2017wasserstein}
A.~Ramdas, N.~Garc{\'\i}a~Trillos, and M.~Cuturi.
\newblock On wasserstein two-sample testing and related families of
  nonparametric tests.
\newblock \emph{Entropy}, 19\penalty0 (2):\penalty0 47, 2017.
\newblock URL \url{https://www.mdpi.com/1099-4300/19/2/47}.

\bibitem[Revay and Teschke(2019)]{revay2019multiclass}
S.~Revay and M.~Teschke.
\newblock Multiclass language identification using deep learning on spectral
  images of audio signals.
\newblock \emph{arXiv preprint arXiv:1905.04348}, 2019.
\newblock URL \url{https://arxiv.org/abs/1905.04348}.

\bibitem[Rote and Zachariasen(2007)]{RoteZachariasen}
G.~Rote and M.~Zachariasen.
\newblock Matrix scaling by network flow.
\newblock In N.~Bansal, K.~Pruhs, and C.~Stein, editors, \emph{Proceedings of
  the Eighteenth Annual {ACM-SIAM} Symposium on Discrete Algorithms, {SODA}
  2007, New Orleans, Louisiana, USA, January 7-9, 2007}, pages 848--854.
  {SIAM}, 2007.
\newblock URL \url{http://dl.acm.org/citation.cfm?id=1283383.1283474}.

\bibitem[Scetbon and Cuturi(2020)]{NEURIPS2020_9bde76f2}
M.~Scetbon and M.~Cuturi.
\newblock Linear time sinkhorn divergences using positive features.
\newblock In H.~Larochelle, M.~Ranzato, R.~Hadsell, M.~Balcan, and H.~Lin,
  editors, \emph{Advances in Neural Information Processing Systems}, volume~33,
  pages 13468--13480. Curran Associates, Inc., 2020.
\newblock URL
  \url{https://proceedings.neurips.cc/paper/2020/file/9bde76f262285bb1eaeb7b40c758b53e-Paper.pdf}.

\bibitem[Scetbon et~al.(2021)Scetbon, Cuturi, and
  Peyr{\'e}]{pmlr-v139-scetbon21a}
M.~Scetbon, M.~Cuturi, and G.~Peyr{\'e}.
\newblock Low-rank sinkhorn factorization.
\newblock In M.~Meila and T.~Zhang, editors, \emph{Proceedings of the 38th
  International Conference on Machine Learning}, volume 139 of
  \emph{Proceedings of Machine Learning Research}, pages 9344--9354. PMLR,
  18--24 Jul 2021.
\newblock URL \url{https://proceedings.mlr.press/v139/scetbon21a.html}.

\bibitem[Schmitzer(2019)]{schmitzer2019stabilized}
B.~Schmitzer.
\newblock Stabilized sparse scaling algorithms for entropy regularized
  transport problems.
\newblock \emph{SIAM Journal on Scientific Computing}, 41\penalty0
  (3):\penalty0 A1443--A1481, 2019.
\newblock \doi{10.1137/16M1106018}.

\bibitem[Schrieber et~al.(2016)Schrieber, Schuhmacher, and
  Gottschlich]{schrieber2016dotmark}
J.~Schrieber, D.~Schuhmacher, and C.~Gottschlich.
\newblock Dotmark--a benchmark for discrete optimal transport.
\newblock \emph{IEEE Access}, 5:\penalty0 271--282, 2016.
\newblock \doi{10.1109/ACCESS.2016.2639065}.

\bibitem[Sinkhorn(1967)]{Sinkhorn1967a}
R.~Sinkhorn.
\newblock Diagonal equivalence to matrices with prescribed row and column sums.
\newblock \emph{The American Mathematical Monthly}, 74\penalty0 (4):\penalty0
  402, 1967.
\newblock \doi{10.2307/2314570}.

\bibitem[Sinkhorn and Knopp(1967)]{Sinkhorn1967}
R.~Sinkhorn and P.~Knopp.
\newblock Concerning nonnegative matrices and doubly stochastic matrices.
\newblock \emph{Pacific Journal of Mathematics}, 21:\penalty0 343--348, 1967.
\newblock ISSN 0030-8730.

\bibitem[Tai et~al.(2021)Tai, Bailis, and Valiant]{pmlr-v139-tai21a}
K.~S. Tai, P.~D. Bailis, and G.~Valiant.
\newblock Sinkhorn label allocation: Semi-supervised classification via
  annealed self-training.
\newblock In M.~Meila and T.~Zhang, editors, \emph{Proceedings of the 38th
  International Conference on Machine Learning}, volume 139 of
  \emph{Proceedings of Machine Learning Research}, pages 10065--10075. PMLR,
  18--24 Jul 2021.
\newblock URL \url{https://proceedings.mlr.press/v139/tai21a.html}.

\bibitem[Villani(2003)]{Villani2003}
C.~Villani.
\newblock \emph{Topics in {O}ptimal {T}ransportation}, volume~58 of
  \emph{Graduate Studies in Mathematics}.
\newblock American Mathematical Society, Providence, RI, 2003.
\newblock ISBN 0-821-83312-X.
\newblock \doi{10.1090/gsm/058}.
\newblock URL \url{http://books.google.com/books?id=GqRXYFxe0l0C}.

\bibitem[Villani(2009)]{Villani2009}
C.~Villani.
\newblock \emph{Optimal transport, old and new}, volume 338 of
  \emph{Grundlehren der Mathematischen Wissenschaften}.
\newblock Springer, Berlin, 2009.
\newblock \doi{10.1007/978-3-540-71050-9}.

\end{thebibliography}

\end{document}